\documentclass[10pt]{article}

\usepackage{amssymb}
\usepackage{amsmath}
\usepackage{amsthm}
\usepackage{mathrsfs}
\usepackage{slashbox}
\usepackage[all]{xy}
\usepackage{graphicx}
\usepackage{wrapfig}

\def\2{{1\over 2}}

\newcommand{\rf}[1]{(\ref{#1})}
\def\b{\bar}

\newcommand{\ud}{\mathrm{d}}
\renewcommand{\t}{\tilde}
\newcommand{\p}{\partial}

\newcommand{\mA}{\mathbf{A} }
\newcommand{\mB}{\mathbf{B} }

\newcommand{\mF}{\mathcal{F} }

\newcommand{\mg}{\mathfrak{g}}

\newcommand{\mG}{\mathfrak{G}}

\title{\bf{Quasiclassical Lian-Zuckerman Homotopy Algebras, Courant Algebroids and Gauge Theory}}
\author{Anton M. Zeitlin\footnote{anton.zeitlin@yale.edu, http://math.yale.edu/$\sim$az84, http://www.ipme.ru/zam.html} \\
Department of Mathematics,\\
Yale University,\\
442 Dunham Lab, 10 Hillhouse Ave\\
New Haven, CT 06511\\} 

\date{}

\begin{document}
\maketitle
\begin{center}
{\it To Gregg J. Zuckerman for his 60th birthday} 
\end{center}
\vspace*{3mm}
 
\begin{abstract}
We define a quasiclassical limit of the Lian-Zuckerman homotopy BV algebra 
(quasiclassical LZ algebra) on the 
subcomplex, corresponding to "light modes", i.e. the elements of zero conformal weight, 
of the semi-infinite (BRST) cohomology complex of the Virasoro algebra associated with vertex operator algebra (VOA) with a formal parameter. We also construct a certain deformation of the BRST differential parametrized by a constant two-component tensor, such that it leads to the deformation of the $A_{\infty}$-subalgebra of the quasiclassical LZ algebra. 
Altogether this gives a functor the category of VOA with a formal parameter to the category of $A_{\infty}$-algebras. 
The associated generalized Maurer-Cartan equation gives the analogue of the Yang-Mills equation for a wide class of VOAs. Applying this construction to an example of VOA generated by $\beta$-$\gamma$ systems, we find 
a remarkable relation between the Courant algebroid and the homotopy algebra of the Yang-Mills theory. 
\end{abstract}
\section{Introduction}
The relation between the dynamics of two-dimensional world and $D$-dimensional field theory is in the very heart of String Theory. An important problem is to find how the classical nonlinear equations of motion of gauge theory and gravity emerge from such two-dimensional dynamics. 

In the early days of string theory, the solutions of the 
linearized equations of motion modulo gauge symmetry were identified with the semi-infinite cohomology classes of the Virasoro algebra for a certain Virasoro module. On the other hand, in the same time period the nonlinear equations of motion, e.g. Yang-Mills and Einstein equations with extra fields, have been derived as relations coming from the conformal invariance condition for two-dimensional sigma models \cite{fts}, \cite{eeq}, \cite{polbook}. These relations are usually written as 
$\beta(\phi_{\{\nu\}},h)=0$, where $\beta(\phi_{\{\nu\}},h)=\sum_{i\in\mathbb{N}}h^i\beta_i(\phi_{\{\nu\}})$ 
is some function of all the fields $\phi_{\nu}$, which are present in the sigma-model, and $h$ is a formal parameter having a meaning of the Planck constant. The equation $\beta_1(\phi_{\{\nu\}})=0$ ($\beta_1$ is known as a 1-loop $\beta$-function) is equivalent to the classical nonlinear field equations. 

The definition of $\beta$ was and is still a mystery from the mathematical point of view. However, it was suggested by some authors in the 1980s that the algebraic version of the equation $\beta(\phi_{\{\nu\}},h)=0$ should be something like the (generalized) Maurer-Cartan equation. 

Soon after that, the open and the closed String Field Theories (SFT) have been constructed \cite{witopen}, \cite{zwiebach}. The equations of motion of String Field Theory, which are supposed to contain nonlinear field-theoretic equations of motion, had the form of Maurer-Cartan equations for certain homotopy algebras: associative algebras for open strings and strong homotopy Lie algebras for closed ones. However, it was very hard to derive these equations from SFT \cite{taylor},\cite{berkovits}, moreover, the relation of these Maurer-Cartan equations to the $\beta$-function vanishing condition was also not clear. It is worth noting several attempts to establish such a relation \cite{banks},\cite{sen}.  
 
At the same time, in the 1980s, the Vertex Operator Algebra (VOA) theory was constructed (see e.g. \cite{fhl}, \cite{benzvi}): a mathematical theory describing two-dimensional conformally  invariant models, 
which are of special importance in string theory. It was observed by B.H. Lian  and G.J. Zuckerman \cite{lz}, 
that a special class of VOAs, the $topological$ VOAs (TVOAs), possesses a  homotopy algebra. 
This algebra turned out to be a homotopy Batalin-Vilkovisky (BV) algebra. 
It was also conjectured in \cite{zuck} that there are "higher homotopies" for this algebra, such that it can be extended to the object which the authors of \cite{zuck} called $G_{\infty}$-algebra (see also \cite{huang}, \cite{gorbounov}), which first appeared in \cite{gj} (see also \cite{tamarkin}). 
In a recent article \cite{vallette}, it was also conjectured and proven (for a certain class of TVOA) that there exists a $BV_{\infty}$-algebra, which is the extension of Lian-Zuckerman homotopy algebra. 

One of the important classes of TVOAs is the semi-infinite cohomology complex (or simply BRST complex)  
associated with certain VOAs. 

Taking into account the physical motivation, 
in order to construct the $\beta$-function and classical field equations,
 one has to construct a functor from the category of vertex operator algebras to the category of homotopy algebras. 
A natural candidate for this functor is the one, provided by the construction of LZ algebra on the BRST complex of the Virasoro algebra. We show by means of several examples, that it is enough to consider a certain quasi-isomorphic subcomplex of this BRST complex, the so-called $light$ $modes$, which are annihilated by the $L_0$ Virasoro mode and therefore generate a subalgebra in the LZ algebra. 

 When conformal weights in the VOA are bounded from below (in this paper we assume them to be bounded by zero), the complex of light modes is easy to work with. In this paper we will consider only this light mode complex. 

 We define what we call the $quasiclassical$ $limit$ of the LZ homotopy BV algebra, which leads to a certain "truncation" of higher homotopies. This means that it contains $A_{\infty}$- and $L_{\infty}$-algebras, such that all polylinear operations vanish starting from the quadrilinear ones. We conjecture that it is actually the $BV_{\infty}$-algebra, motivated by the results of \cite{gorbounov}, \cite{vallette}. In the $\beta$-$\gamma$ example, this gives a homotopy BV algebra associated to Courant algebroid on the sum of tangent and cotangent bundles. The corresponding $L_{\infty}$-algebra is the one, constructed by D. Roytenberg and A. Weinstein \cite{weinstein}. 

Then we introduce a deformation of the BRST differential, 
which we call a $flat$ $background$ deformation, which 
corresponds to the Abelian vertex subalgebra and involves a constant two-component tensor. 
We prove that this deformation can be continued to the deformation of the homotopy commutative $A_{\infty}$ 
subalgebra of the LZ algebra. This gives us a functor from the category of VOAs with a chosen Abelian subalgebra 
to the category of $A_{\infty}$-algebras parametrized by this constant tensor. In the beta-gamma example, this deformed $A_{\infty}$-algebra turns out to be the $A_{\infty}$-algebra of the Yang-Mills theory with external fields.  

The structure of the paper is as follows. In Section 2 we set up notations and recall the basic facts concerning  Lian-Zuckerman homotopy BV algebra. We also give a short reminder about $A_{\infty}$-algebras.  In Section 3, we define the quasiclassical  limit of the LZ algebra of light modes. We prove explicitly that the quasiclassical limit of LZ product  extends to the $A_{\infty}$-algebra. In Section 4 we define what we call a {\it flat background} deformation of the quasiclassical LZ algebra. This is a deformation of the 
differential, such that in general only the homotopy commutative  $A_{\infty}$-subalgebra of the homotopy 
BV algebra can be deformed in such a way to satisfy all necessary conditions with the deformed differential. 

In Section 5 we apply the constructions we introduced earlier to the VOA generated by a family of $\beta$-$\gamma$ systems. In particular, we obtain that the quasiclassical LZ algebra corresponds to the BV algebra, which in conformal weight 1 reproduces  the Courant algebroid.  After that, we show that the flat background deformation in a certain case leads to the Yang-Mills theory with matter fields,  which yields a remarkable relation between Courant/Dorfman brackets and gauge theory. In the last section we outline some of the numerous possible directions continuing the studies started in this paper.

\section{Lian-Zuckerman homotopy BV algebra.}
\noindent{\bf 2.1. Notation and Conventions.} Throughout the paper we will work with vertex operator algebras 
(VOA), using physics notation. Therefore  the elements of the VOA's vector space will be referred to as $states$ and $A(z)$ denotes the vertex operator $Y(A,z)$ (see e.g. \cite{benzvi}, \cite{fhl}), corresponding to the state $A$. To simplify the calculations, we introduce a special notation for certain operator product coefficients. Namely if $A$ and $B$ are the elements of the VOA, then we denote  
\begin{eqnarray}
\langle A,B\rangle\equiv Res_z(zA(z)B),\quad [A,B]\equiv Res_z(A(z)B), \quad AB\equiv Res_{z}\bigg(\frac{A(z)B}{z}\bigg).
\end{eqnarray}

\noindent{\bf 2.2. Topological VOA and the Lian-Zuckerman homotopy BV algebra.} Topological vertex operator algebra (TVOA) is a vertex superalgebra (see e.g.\cite{benzvi}) that has an additional odd operator 
$Q$ which makes the graded vector space of VOA a chain complex, such that the Virasoro element $L(z)$ is $Q$-exact. The formal definition is as follows (see e.g. \cite{zuck} for more details).\\

\noindent{\bf Definition 2.1.} {\it Let V be a $\mathbb{Z}$-graded vertex operator superalgebra, such that $V=\oplus_i V^i=\oplus_{i,\mu}V^i[\mu]$, where $i$ represents grading of $V$ with respect to conformal weight and $\mu$ 
represents fermionic grading of $V^i$. We call V a topological vertex operator algebra (TVOA) if there exist four elements: $J\in V^1[1]$, $b\in V^2[-1]$, $F\in V^1[0]$, $L\in V^2[0]$, such that 
\begin{eqnarray}
[Q,G(z)]=\mathcal{L}(z),\quad  Q^2=0,\quad G_0^2=0,
\end{eqnarray}
where $Q=J_0$ and $G(z)=\sum_n b_nz^{-n-2}$, $J(z)=\sum_n J_nz^{-n-1}$, \\
$\mathcal{L}(z)=\sum_n\mathcal{L}_nz^{-n-2}$, 
$F(z)=\sum_nF_nz^{-n-1}$. 
Here $\mathcal{L}(z)$ is the Virasoro element of $V$; the operators $F_0$, $\mathcal{L}_0$ are diagonalizable, commute with each other and 
their egenvalues coincide with fermionic grading and conformal weight correspondingly.}\\

A natural example of such object, which will be crucial in the following, is the semi-infinite cohomology (or simply BRST) complex for the Virasoro algebra \cite{fgz} of some VOA with central charge equal to 26. 
The necessary setup for the construction of the semi-infinite complex (for more details, see \cite{fgz} or \cite{linshaw}) 
is the VOA $\Lambda$,obtained from the following super Heisenberg algebra:
\begin{equation}
\{b_n, c_m\}=\delta_{n+m,0}, \quad n,m \in \mathbb{Z}.
\end{equation}
One can construct the space of $\Lambda$ as a Fock module:
\begin{eqnarray}
\Lambda&=&\{b_{-n_1}\dots b_{-n_k}c_{-m_1}\dots c_{-m_l}\mathbf{1},  n_1,\dots,  n_k>0, m_1,\dots, m_l>0;\nonumber\\
&& c_k \mathbf{1}=0,\ k\geqslant 2;\quad b_k\mathbf{1}=0, \ k\geqslant -1\}.
\end{eqnarray}
Then one can define two fields:
\begin{equation}
b(z)=\sum_m b_mz^{-m-2}, \qquad c(z)=\sum_n c_nz^{-n+1},
\end{equation}
which according to the commutation relations between modes have the following operator product:
\begin{equation}
b(z)c(w)\sim\frac{1}{z-w}.
\end{equation}
The Virasoro element is given by the following expression:
\begin{equation}
L^{\Lambda}(z)=2:\partial b(z)c(z):+:b(z)\partial c(z):,
\end{equation}
such that $b(z)$ has conformal weight $2$, and $c(z)$ has conformal weight $-1$. Here, as usual, $: :$ stand for normal ordered product, e.g. $b(z)c(w)=\frac{1}{z-w}+:b(z)c(w):$ (for more details see e.g. \cite{benzvi}, Section 2.2). 

Now let $V$ be a VOA with the Virasoro element $L(z)$. Let us consider the tensor product $V\otimes\Lambda$. Then we have the following Proposition. 
\\

\noindent {\bf Proposition 2.1.} \cite{fgz} {\it If V is a VOA with the central charge of Virasoro algebra equal to 26, then $V\otimes \Lambda$ is a topological vertex algebra, where 
\begin{eqnarray}
&&J(z)=c(z)L(z)+:c(z)\p c(z) b(z):+\frac{3}{2}\p^2 c(z), \quad G(z)=b(z),\nonumber\\ 
&&F(z)=:c(z)b(z):, \quad \mathcal{L}(z)=L(z)+L^{\Lambda}(z).
\end{eqnarray}}
The operator $Q=J_0$ is traditionally called the $BRST$ $operator$ and the eigenvalue of $F_0$,
 i.e. fermionic grading is usually called the $ghost$ $number$.

Lian and Zuckerman observed that each TVOA possesses a rich algebraic structure. They have shown the following. One can define two operations which are cochain maps with respect to $Q$: 
\begin{eqnarray}\label{mu}
\mu(a_1,a_2)=Res_z\frac{ a_1(z)a_2}{z},\quad \{a_1,a_2\}=\frac{(-1)^{|a_1|}}{2\pi i} \oint dz (b_{-1}a_1)(z)a_2.
\end{eqnarray}
These operations satisfy the following relations.\\

\noindent {\bf Proposition 2.2.}\cite{lz} {\it The operation $\mu$ is homotopy commutative and homotopy associative:

\begin{eqnarray}\label{lzrel}
&&Q\mu(a_1,a_2)=\mu(Q a_1,a_2)+(-1)^{|a_1|}\mu(a_1,Q a_2),\nonumber\\
&&\mu(a_1,a_2)-(-1)^{|a_1||a_2|}\mu(a_2,a_1)= \nonumber\\
&&Qm(a_1,a_2)+m(Qa_1,a_2)+(-1)^{|a_1|}m(a_1,Qa_2),\nonumber\\
&& Qn(a_1,a_2,a_3)+n(Qa_1,a_2,a_3)+(-1)^{|a_1|}n(a_1,Qa_2,a_3)+\nonumber\\
&&(-1)^{|a_1|+|a_2|}n(a_1,a_2,Qa_3)=\mu(\mu(a_1,a_2),a_3)-\mu(a_1,\mu(a_2,a_3))
\end{eqnarray}
where 
\begin{eqnarray}
m(a_1,a_2)=\sum_{i\ge 0}\frac{(-1)^i}{i+1}Res_wRes_{z-w}(z-w)^iw^{-i-1}b_{-1}
(a_1(z-w)a_2)(w)\mathbf{1},
\nonumber\\
n(a_1,a_2,a_3)=\sum_{i\ge 0}\frac{1}{i+1}Res_zRes_w w^iz^{-i-1}(b_{-1}a_1)(z)a_2(w)a_3+\nonumber\\
\ \ \ \ \ \ (-1)^{|a_1||a_2|}\sum_{i\ge 0}\frac{1}{i+1}Res_wRes_z z^iw^{-i-1}(b_{-1}a_2)(w)a_1(z)a_3.
\end{eqnarray}}

\noindent The operation $\{\cdot,\cdot\}$ measures the failure of $b_0$ to be a derivation of $\mu$. In other words we have the following Proposition.\\

\noindent {\bf Proposition 2.3.}\cite{lz} {\it The operations $\mu$ and $\{\cdot,\cdot\}$ are related in the following way:
\begin{eqnarray}\label{lzrel2}
\{a_1,a_2\}=b_0\mu(a_1,a_2)-\mu(b_0a_1,a_2)-(-1)^{|a_1|}\mu(a_1,b_0a_2).
\end{eqnarray}}
\noindent Moreover, the bracket satisfies the relations of a homotopy Gerstenhaber algebra.\\

\noindent {\bf Proposition 2.4.}\cite{lz} {\it The operations $\mu$ and $\{\cdot,\cdot\}$ satisfy the relations: 
\begin{eqnarray}\label{lzrel3}
&&\{a_1,a_2\}+(-1)^{(|a_1|-1)(|a_2|-1)}\{a_2,a_1\}=\\
&&(-1)^{|a_1|-1}(Qm'(a_1,a_2)-m'(Qa_1,a_2)-(-1)^{|a_2|}m'(a_1,Qa_2)),
\nonumber\\
&& \{a_1,\mu(a_2,a_3)\}=\mu(\{a_1,a_2\},a_3)+(-1)^{(|a_1|-1)||a_2|}\mu(a_2,\{a_1, a_3\}),\nonumber\\
&&\{\mu(a_1,a_2),a_3\}-\mu(a_1,\{a_2,a_3\})-(-1)^{(|a_3|-1)|a_2|}\mu(\{a_1,a_3\},a_2)=\nonumber\\
&&(-1)^{|a_1|+|a_2|-1}(Qn'(a_1,a_2,a_3)-n'(Qa_1,a_2,a_3)-\nonumber\\
&&(-1)^{|a_1|}n'(a_1,Qa_2,a_3)-(-1)^{|a_1|+|a_2|}n'(a_1,a_2,Qa_3),\nonumber\\
&&\{\{a_1,a_2\},a_3\}-\{a_1,\{a_2,a_3\}\}+(-1)^{(|a_1|-1)(|a_2|-1)}\{a_2,\{a_1,a_3\}\}=0.\nonumber
\end{eqnarray}
where $m', n'$ are some bilinear and trilinear operations on the TVOA (see p. 621 of \cite{lz}).}\\

\noindent In fact, the operations $m'$ and  $n'$ are constructed from $\mu,m$ and $n$. However, we will not need explicit expressions in the following. The relations \rf{lzrel}-\rf{lzrel3} by definition mean that $\mu(\cdot,\cdot)$ and $\{\cdot,\cdot \}$ generate the $homotopy$ $BV$ $algebra$.  In the following we will refer to the concrete homotopy BV algebra generated by $\mu$ and $\{,\}$ as $LZ$  $algebra$. 

We note here, that in this article when we say "homotopy" (associative, Lie, Gerstenhaber, BV...) algebra it means that we have certain set of operations endowed with just a first level of homotopies. This provides the structure of a strict (associative, Lie, Gerstenhaber,  BV...) algebra structure on the cohomology, as in \cite{lz}. When we talk about 
$\infty$-algebras, we mean that there is a full set of higher homotopies (see below the explicit description of 
$A_{\infty}$-algebras) satisfying some relations. Therefore $\infty$-algebra is a homotopy algebra, but the inverse is not true in general.

Lian and Zuckerman made a conjecture, that the LZ algebra can be extended to $G_{\infty}$-algebra (see e.g. \cite{zuck}). It was recently proved \cite{gorbounov}, \cite{vallette} that for a certain class of TVOAs, that this algebra is in fact the $BV_{\infty}$-algebra \cite{vallette}. As we discussed above, this means that there are "higher homotopies", i.e. in general case nonzero multilinear operations which satisfy the higher associativity/Jacobi/Leibniz relations.

Here we will not discuss this complicated object in detail. However, this conjecture implies that there 
exists a homotopy algebra, which "extends" the relations between $\mu$ and $n$ only. Such an algebra is called 
an $A_{\infty}$-algebra, and we will discuss the precise definition in the next subsection. \\

\noindent{\bf 2.3. Short reminder of $A_{\infty}$-algebras.}
The $A_{\infty}$-algebra is a generalization of a differential graded associative algebra. Namely, consider a graded vector space $V$ with a differential $Q$. Consider the multilinear operations $\mu_i: V^{\otimes i}\to V$ of the degree $2-i$, such that $\mu_1=Q$. \\

\noindent {\bf Definition 2.2.} (see e.g. \cite{stashbook}){\it The space V is an $A_{\infty}$-algebra if the operations 
$\mu_n$ satisfy bilinear identity:
\begin{eqnarray}\label{arel}
\sum^{n-1}_{i=1}(-1)^{i}M_i\circ M_{n-i+1}=0 
\end{eqnarray}
on $V^{\otimes n}$,  
where $M_s$ acts on $V^{\otimes m}$ for any $m\ge s$ as the sum of all possible operators of the form 
${\bf 1}^{\otimes^l}\otimes\mu_s\otimes{\bf 1}^{\otimes^{m-s-l}}$ taken with appropriate signs. In other words, 
\begin{eqnarray}
M_s=\sum^{n-s}_{l=0}(-1)^{l(s+1)}{\bf 1}^{\otimes^l}\otimes\mu_s\otimes{\bf 1}^{\otimes^{m-s-l}}.
\end{eqnarray}}
Let us write several relations which are satisfied by $Q$, $\mu_1$, $\mu_2$, $\mu_3$:
\begin{eqnarray}
&&Q^2=0,\\
&&Q\mu_2(a_1,a_2)=\mu_2(Q a_1,a_2)+(-1)^{|a_1|}\mu_2(a_1,Q a_2),\nonumber\\
&&Q\mu_3(a_1,a_2, a_3)+\mu_3(Q a_1,a_2, a_3)+(-1)^{|a_1|}\mu_3(a_1,Q a_2, a_3)+\nonumber\\
&&(-1)^{|a_1|+|a_2|}\mu_3( a_1, a_2, Q a_3)=\mu_2(\mu_2(a_1,a_2),a_3)-\mu_2(a_1,\mu_2(a_2,a_3)).\nonumber
\end{eqnarray}

In such a way we see that if $\mu_n=0$, $n\ge 3$ , then we have just a differential graded associative algebra (DGA).
If the operations $\mu_n$ vanish for all $n>k$, such $A_{\infty}$-algebras are sometimes called 
$A_{(k)}$-algebras \cite{stasheff}, so e.g. DGA is $A_{(2)}$ algebra.   

We observe that putting $\mu_2\equiv \mu$ and $\mu_3=n$, these relations are manifestly the same as 
the ones relating $Q$, $\mu$ and $n$. The Lian-Zuckerman conjecture states that there are "higher homotopies" $\mu_n, n>3$ satisfying the relations \rf{arel}.   

It is well known that the relations \rf{arel} can be encoded into one equation $\p^2=0$ \cite{stashbook}. To see this one can apply the desuspension operation (the operation which shifts the grading $s^{-1}: V_{n}\to (s^{-1}V)_{n-1}$) to $\mu_n$. In such a way we can define operations of degree $1$: $\t \mu_n=s\mu_n {(s^{-1})}^{\otimes n}$. More explicitly, 
\begin{eqnarray}
\t \mu_n(s^{-1} a_1,...,s^{-1} a_n)=(-1)^{s(a)}s^{-1}\mu_n(a_1,...,a_n),
\end{eqnarray}
such that $s(a)=(1-n)|a_1|+(2-n)|a_2|+...+|a_{n-1}|$. 
The relations between $\t\mu_n$ operations can be summarized in the following simple equations:
\begin{eqnarray}\label{M2}
\sum^n_{i=1}\t M_i\circ \t M_{n+1-i}=0.
\end{eqnarray}
on $V^{\otimes n}$, where each $\t M_s$ acts on $V^{\otimes m}$ (for $m\ge s$) as the sum of all operators ${\bf 1}^{\otimes l}\otimes\t \mu_s\otimes{\bf 1}^{\otimes k}$, such that $l+s+k=m$. 
Combining them into one operator $\p=\sum_n\t M_n$, acting on a space $\oplus_kV^{\otimes k}$ the relations \rf{arel} can be summarized in one equation $\p^2=0$. 

An important object in the theory of $A_{\infty}$-algebras is the generalized Maurer-Cartan (GMC) equation. Let us pick 
$X\in V$ of degree 1. 
Then the equation 
\begin{eqnarray}
QX+\sum_{n\ge2}\mu_n(X,...,X)=0
\end{eqnarray}
is called the \emph{generalized Maurer-Cartan equation} on $X$. It is worth mentioning that it is well defined in general only on nilpotent elements, i.e. such that $\mu_n(X,\dots, X)=0$ for $n>k$. This will not be a problem in the following, because the only $A_{\infty}$ algebras we will consider in this article, are $A_{(3)}$-algebras and therefore GMC equation will be well defined for all elements of degree 1. 
 
Generalized Maurer-Cartan equation is known to have the following infinitesimal symmetry:
\begin{eqnarray}
X\mapsto X +\epsilon (Q\alpha+\sum_{n\ge2,k}(-1)^{n-k}\mu_n(X,...,\alpha,...,X)), 
\end{eqnarray}
where $\epsilon$ is infinitesimal, $\alpha$ is an element of degree 0 and $k$ means the position of $\alpha$ in $\mu_n$. \\
 
\section{Light modes and the quasiclassical limit}
\noindent{\bf 3.1. Vertex operator algebras with a formal parameter.}
In this article we will be interested in computing quasiclassical limits, therefore we need to consider the vertex operator algebras depending on a formal parameter. Namely, we consider the VOAs on the spaces of the form $\mathcal{V}=\oplus_{n\in \mathbb{Z}_{\ge 0}}\mathcal{V}^n$, where $n$ stands for the grading with respect to conformal weight, 
and $\mathcal{V}^n=V^n[h^{-1},h]$, where $V^n$ are some vector spaces and $h$ is a formal parameter. Let us denote $V=\oplus_{n\in \mathbb{Z}_{\ge 0}}V^n$. 

\noindent Then, we require the operator products to meet the following conditions: \\

\noindent{\it i)  For any state $A$ the associated vertex operator $A(z)=\sum_nA_nz^{-n-1}$ is such that its Fourier modes $A_n\in End_{\mathbb{C}[h,h^{-1}]}(\mathcal{V})$, i.e. they commute with the natural action of $\mathbb{C}[h,h^{-1}]$ on $\mathcal{V}$.\\
  
\noindent ii) Let $A,B\in V$ and $A(z)=\sum_nA_nz^{-n-1}$,  then
\begin{eqnarray}
A_nB\in hV[h],  \quad n\ge 0. 
\end{eqnarray}}
Moreover, we put the following conditions on the Virasoro action:\\

\noindent {\it iii)  $V$ is invariant under the action of the operator $L_{-1}$. \\

\noindent iv) Let $A\in V$ be a state of conformal weight $1$, then $L_1A\in hV$.}
\footnote[2]{This condition, in principle, can be relaxed, by letting $L_1A\in hV[h]$, and (with some modifications) one can apply the constructions of the article to this case. However, being motivated by certain concrete examples we keep the more restrictive form of condition $iv)$.}\\

\noindent As a consequence of properties $i)-iv)$, we observe that letting $h=1$ we obtain a VOA structure on $V$. 

\noindent Let us consider some simple examples of such parameter-dependent VOAs.\\

\noindent{\bf Examples.}\\ 

\noindent {\it i) Heisenberg VOA.} We consider the Fock space for the Heisenberg algebra, 
\begin{eqnarray}
[a_n,a_m]=hn\delta_{n,-m}
\end{eqnarray}
i.e. the space $F_{a}[h^{-1},h]$, where $F_{a}=\{a_{-n_1}...a_{-n_k}|0\rangle, n_1, ..., n_k>0 ; a_n|0\rangle=0, n\ge 0\}$.  
The space  $F_{a}$ is a VOA, such that on the operator product expansion language these 
commutation relations are summarized as follows:
\begin{eqnarray}
a(z)a(w)\sim \frac{h}{(z-w)^2},
\end{eqnarray}
where $a(z)=\sum_na_nz^{-n-1}$. The Virasoro element is given by the operator $L(z)=\frac{1}{h}:a(z)^2$:, where dots stand for standard normal ordering in the Fock space.\\

\noindent 
{\it ii) $\beta$-$\gamma$ system.} This is the most interesting example for us, because in this case there is  
a nontrivial subspace of the VOA of conformal weight equal to zero. 
For all the details we refer to the paper $\cite{cdr}$. 
Let us consider the Heisenberg algebra of the following kind:
\begin{eqnarray}\label{heis}
[\gamma_n,\beta_m]=h\delta_{n,-m}.
\end{eqnarray}
The construction of the corresponding VOA space $F_{\beta\gamma}[h^{-1},h]$ is as follows.   
Let us consider again the Fock space for the Heisenberg algebra \rf{heis}:  
$\t F_{\beta\gamma}=\{\beta_{-n_1}....\beta_{-n_k}\gamma_{-m_1}...\gamma_{-m_l}|0\rangle, n_1,..., n_k>0, m_1, ..., m_l>0;$ 
$\beta_n|0\rangle=0, n\ge 0, \gamma_n |0\rangle=0, n>0\}$. In other words, $\t F_{\beta\gamma}=
\mathbb{C}[\beta_{-n}, \gamma_{-m}]_{n>0, m\ge 0}$. The space  $\t F_{\beta\gamma}[h^{-1},h]$ already carries an algebraic structure of a VOA, but one can proceed further and construct the space 
$F_{\beta\gamma}=\t F_{\beta\gamma}\otimes_{A} \hat{A}$, where 
$A=\mathbb{C}[\gamma_0]$ and $\hat{A}=\mathbb{C}[[\gamma_0]]$ or any other function field containing 
$\mathbb{C}[\gamma_0]$. From \cite{cdr} one can prove that $F_{\beta\gamma}[h^{-1},h]$ is a VOA with a formal parameter, such that the commutation relations \rf{heis} can be expressed via the operator product:
\begin{eqnarray}
\gamma(z)\beta(w)\sim\frac{h}{z-w},\quad  \beta(z)\beta(w)\sim 0, \quad \gamma(z)\gamma(w)\sim 0,
\end{eqnarray}
where $\beta(z)=\sum_n\beta_nz^{-n-1}$ is a quantum field of conformal dimension $1$ and 
$\gamma(z)=\sum_n\gamma_nz^{-n}$ is a quantum field of conformal dimension $0$. The Virasoro element is given by 
$L_{\beta,\gamma}(z)=-\frac{1}{h}:\beta\p \gamma(z):$, where dots stand for normal ordering. 
\\

One can easily continue this list of examples. In fact, there is a quite large class of VOAs, 
which can be obtained from the VOAs with a formal parameter by letting this parameter be equal to $1$.\\

\noindent{\bf 3.2. The light modes.} Let us consider the space 
$\mathcal{V}\otimes \Lambda$, where $\Lambda$ is the $b$-$c$ ghost VOA defined in subsection 2.2. When the central charge of the Virasoro algebra of $\mathcal{V}$ is equal to 26, this space is the semi-infinite cohomology complex of the Virasoro algebra. However, there is a subspace of   $\mathcal{V}\otimes \Lambda$ which is a complex with respect to BRST operator regardless of the value of the central charge.\\

\noindent {\bf Definition 3.1.}{\it We call the subspace of $\mathcal{V}\otimes \Lambda$ which is annihilated by $\mathcal{L}_0$  (i.e. space of states of conformal weight zero) the space of {\rm light modes}.} \\

\noindent Since the conformal weights in the VOA are greater than zero, the following proposition holds.\\

\noindent{\bf Proposition 3.1.}{\it (i)The space of light modes $F_{\mathcal{L}_0}$ 
is linearly spanned by the elements which correspond to the operators:
\begin{eqnarray}
&&u(z), \quad c(z)A(z), \quad \p c(z) a(z), \quad c(z)\p c(z) \t A(z), \nonumber\\
&& c(z)\p^2 c(z) \t a(z), \quad c(z)\p c(z)\p^2 c(z) \t u(z).  
\end{eqnarray}
Here $u, \t u$, $a, \t a\in \mathcal{V}$ are of conformal weight 0 and $A, \t A\in\mathcal{V}$ are of conformal weight 1. \\
(ii)The space $F_{\mathcal{L}_0}$ is a chain complex, quasi-isomorphic to the semi-infinite complex when the central charge is equal to 26. The differential acts on $F_{\mathcal{L}_0}$ in the 
following way (we recall that the grading is given by the ghost number):

\begin{eqnarray}
\xymatrixcolsep{30pt}
\xymatrixrowsep{3pt}
\xymatrix{
& \mathcal{V}^{1}\ar[ddddr] & \mathcal{V}^{1}\ar[ddddr]^{\frac{1}{2}L_{1}}& \\
&& L_{-1} &&\\
& \bigoplus & \bigoplus & \\
&& {-\frac{1}{2}L_{1}} &&\\
\mathcal{V}^{0}\ar[uuuur]^{L_{-1}} & \mathcal{V}^{0}\ar[uuuur]\ar[r]_{i\rm{d}} & \mathcal{V}^{0}  & \mathcal{V}^{0}
}
\end{eqnarray}
where $\mathcal{V}^{i}$ $(i=0,1)$ is the space of the elements of $\mathcal{V}$ of conformal dimension $i$. 
}\\

\noindent{\bf Proof.} The first part of the Proposition can be proven by comparing the conformal weights of differential polynomials of operator $c(z)$ and the operators in vertex algebra $V$ in such a way that the total conformal weight is zero. Part ii) can be obtained using the following formulas:
\begin{eqnarray}\label{addf}
QA(z)=\p c(z) A(z)+c(z)\p A(z) +\frac{1}{2}\p^2c(z)L_1A(z), \quad Qc(z)=c(z)\p c(z). 
\end{eqnarray}
where $u\in V^0$, $A\in V^1$. Really, using \rf{addf}, we get:
\begin{eqnarray}
&&Qu(z)=c(z)\p u(z), \quad Q(c(z)A(z))=-c(z)\p^2c(z)L_1A(z),\\ 
&&Q(\p c(z)a(z))=c(z)\p^2 c(z)a(z)+c(z)\p c(z)\p a(z),\nonumber\\ 
&&Q(c(z)\p c(z)\t A(z))=\frac{1}{2}c(z)\p c(z)\p^2 c(z)L_1\t A(z), \quad Q(c(z)\p^2c(z)\t a(z))=0\nonumber, 
\end{eqnarray}
The the Proposition 3.1 is proved.   $\hfill \blacksquare$\\

\noindent Another important observation is about operator $b_0$.\\

\noindent{\bf Corollary 3.1.} 
{\it The space $F_{\mathcal{L}_0}$ is invariant under the action of the operator $b_0$, 
$[Q,b_0]=0$ on $F_{\mathcal{L}_0}$ and $b_0$ has a trivial cohomology. The explicit form of the action of $b_0$ is:
\begin{eqnarray}
\xymatrixcolsep{30pt}
\xymatrixrowsep{3pt}
\xymatrix{
& \mathcal{V}^{1} & \mathcal{V}^{1}\ar[l]_{-i\rm{d}}& \\
& \bigoplus & \bigoplus & \\
\mathcal{V}^{0} & \mathcal{V}^{0}\ar[l]_{i\rm{d}} & \mathcal{V}^{0}  & \mathcal{V}^{0}\ar[l]_{-i\rm{d}}
} 
\end{eqnarray}}

\noindent {\bf Important remark about notation.} In order to simplify some of the calculations, in the following 
we will sometimes write instead of the element of the complex $F_{\mathcal{L}_0}$ the corresponding state of $\mathcal{V}$, i.e. instead of the element corresponding to $c A$ we will just write $A$, or $\t A$ instead of $c\p c\t A$.  Since according to our notation the states without tilde correspond to tensor products with elements of ghost number 0 and 1, and the states marked by tilde correspond to the tensor product with elements of ghost number 2 and 3, this should never lead to confusion.\\

\noindent 
We note here that the operations $\mu$ and $\{\cdot,\cdot\}$ act invariantly on the space of light modes, therefore the light modes form a homotopy subalgebra of LZ algebra. Moreover, from Proposition 3.1. and Corollary 3.1. we see that the Lian-Zuckerman algebra on light modes is determined by means of operator product expansion of the elements of $\mathcal{V}$ of conformal dimensions 0 and 1.  
For example, explicit expressions for the bilinear operation $\mu$ is collected in the table:

\begin{equation}\label{mutab}
\mu(a_1,a_2)=
\end{equation}
\begin{equation}
\begin{tabular}{|l|c|c|c|c|c|r|}
\hline
 \backslashbox{$a_2$}{$a_1$}&          $u_1$ & $A_1$ & $v_1$ &$\t A_1$ &$\t v_1$& $\t u_1$ \\
\hline
$u_2$ &                  $u_1u_2$    &$A_1u_2$ &$v_1u_2$&
$\t A_1u_2$&$\t v_1 u_2$ &$\t u_1 u_2$\\
&  & $+[A_1,u_2]$  &  &   && \\
\hline
$A_2$     &  $u_1A_2$ & $-[A_1,A_2]+$ &$-v_1A_2$  
&  $\frac{1}{2}\langle\t A_1,A_2\rangle$&$[\t v_1,A_2]$&0\\                    
& & $\frac{1}{2}\langle A_1, A_2\rangle$  & $-[\t v_1,A_2]$  &   && \\
\hline
$v_2$& $u_1\t u_2$ &  $A_1v_2$ & 0 & $-[\t A_1, v_2] $ &$-\t v_1v_2$& 0\\
\hline
$\t A_2$ & $u_1\t A_2$ &$\frac{1}{2}\langle A_1,\t A_2\rangle$&  
$[v_1,\t A_2]$ & 0& 0&0\\
\hline
$\t v_2$& $u_1\t u_2$ &   $[A_1,\t v_2]$& $-v_1\t v_2$ & 0  &0& 0\\
\hline
$\t u_2$    & $u_1\t u_2$ &   0& 0 & 0  &0& 0\\
\hline
\end{tabular}\nonumber\\
\end{equation}
\vspace{3mm}

Let us denote by $F^+_{\mathcal{L}_0}$ the space of light modes, which belong to $V[h]\otimes \Lambda$. 
The Lian-Zuckerman algebraic operations act on this space invariantly since $V[h]$ is a vertex algebra. Let us denote $F^+_{\mathcal{L}_0}(1)$ the Lian-Zuckerman algebra on $F^+_{\mathcal{L}_0}$ when $h=1$. Below we construct the embedding of  $F^+_{\mathcal{L}_0}(1)$ into $F^+_{\mathcal{L}_0}$ as a chain complex. 
The following Proposition holds.\\

\noindent{\bf Proposition 3.2.}
{\it Let $\mF$ be the subspace of $F^+_{\mathcal{L}_0}$, linearly 
spanned by the elements corresponding to the vertex operators: 
\begin{eqnarray}
&&u(z),\quad c(z)A(z),\quad h\p c(z) v(z), \quad hc\p c(z) \t A(z),\nonumber\\
&& hc(z)\p^2 c(z) \t v(z),\quad h^2c(z)\p c(z)\p^2 c(z) \t u(z),
\end{eqnarray}
where $u,\t u, v, \t v , A, \t A\in V$. Then the following two statements hold:\\
i) The complex $\mF$ is a BRST subcomplex of $F^+_{\mathcal{L}_0}$ isomorphic to $F^+_{\mathcal{L}_0}(1)$. 
Moreover, it is a subcomplex with respect to the operator $h^{-1}b_0$.\\
ii) The Lian-Zuckerman operations act as follows on $\mF$:
\begin{eqnarray}
&&\mu: \mF_i\otimes \mF_j\to \mF_{i+j}[h], \quad  m: \mF_i\otimes \mF_j\to \mF_{i+j-1}[h], \nonumber\\
&& n: \mF_i\otimes\mF_j\otimes \mF_k\to \mF_{i+j+k-1}[h], 
\quad  \{\mF_i,\mF_j\}\to h\mF_{i+j-1}[h],
\end{eqnarray}
where $\mF_i$ is a subspace of $\mF$ of ghost number $i$.} \\

\noindent Therefore, one can consider the expansions $\mu(\cdot,\cdot)=\mu_0(\cdot,\cdot)+O(h)$, 
$m(\cdot,\cdot)_h=m_0(\cdot,\cdot)+ O(h)$, $n(\cdot,\cdot,\cdot)_h=n_0(\cdot,\cdot,\cdot)+ O(h)$, 
$\{\cdot,\cdot\}=h\{\cdot,\cdot\}_0+O(h^2)$. Then we have a theorem, which follows from Proposition 3.2.\\

\noindent{\bf Theorem 3.1.}{\it The operations $\mu_0(\cdot,\cdot)$, $m_0(\cdot,\cdot)$, $\{\cdot,\cdot\}_0$, $n_0(\cdot,\cdot,\cdot)$ and the differentials $Q$ and $h^{-1}b_0$, defined on the space of the light modes of V ( i.e. 
$F_{\mathcal{L}_0}(1)$), satisfy the relations of the homotopy BV algebra.}\\

\noindent We will call the resulting algebra on $\mF$ the {\it quasiclassical } Lian-Zuckerman algebra \footnote{Let us stress here, that the quasiclassical limit we consider, is different from the standard limit that takes vertex algebra to Poisson vertex algebra. If one applies such limit to our construction, one gets only the subalgebra of the homotopy Lie algebra, which is a part of the homotopy BV algebra.}. 

Further we need explicit expressions for all bilinear operations and their homotopies. 
We use the following notation. Let $\langle a_1,a_2\rangle_0=lim_{h\to 0}h^{-1}\langle a_1,a_2\rangle$ and 
$[ a_1,a_2]_0=lim_{h\to 0}h^{-1}[a_1,a_2]$, where $a_1,a_2\in V$. Below, the expression $a_1 a_2$ means  the normal ordered product of the two elements of $V$ in the $h\to 0$ limit.

Then we can express the bilinear operations $\mu_0(\cdot, \cdot)$ and $\{\cdot,\cdot\}_0$ via the following tables:
\begin{center}
$\mu_0(a_1,a_2)$=
\end{center}
\begin{tabular}{|l|c|c|c|c|c|r|}
\hline
\backslashbox{$ a_2$}{$a_1$}&          $u_1$ & $A_1$ & $v_1$ &$\t A_1$ &$\t v_1$& $\t u_1$ \\
\hline
$u_2$ &                  $u_1u_2$    &$A_1u_2$ &$v_1u_2$&
$\t A_1u_2$&$\t v_1 u_2$ &$\t u_1 u_2$\\
&  & $+[A_1,u_2]_0$  &  &   && \\
\hline
$A_2$     &  $u_1A_2$ & $-[A_1,A_2]_0-$ &$-v_1A_2$  
&  $\frac{1}{2}\langle\t A_1,A_2\rangle_0$&$[\t v_1,A_2]_0$&0\\                    
& & $\frac{1}{2}\langle A_1, A_2\rangle_0$  &   &   && \\
\hline
$v_2$& $u_1\t u_2$ &  $A_1v_2$ & 0 & 0  &$-\t v_1v_2$& 0\\
\hline
$\t A_2$ & $u_1\t A_2$ &$\frac{1}{2}\langle A_1,\t A_2\rangle_0$&  
0 & 0& 0&0\\
\hline
$\t v_2$& $u_1\t u_2$ &   $[A_1,\t v_2]_0$& $-v_1\t v_2$ & 0  &0& 0\\
\hline
$\t u_2$    & $u_1\t u_2$ &   0& 0 & 0  &0& 0\\
\hline
\end{tabular}\\
\vspace{3mm}

\begin{center}
$\{a_1,a_2\}_0$=
\end{center}
\begin{tabular}{|l|c|c|c|c|c|r|}
\hline
\backslashbox{$ a_2$}{$a_1$}&          $u_1$ & $A_1$ & $v_1$ &$\t A_1$ &$\t v_1$& $\t u_1$ \\
\hline
$u_2$ &    0    &$-[A_1,u_2]_0$ & 0 & $[\t A_1,u_2]_0$& 0 & 0\\
\hline
$A_2$     &  0 & $-[A_1,A_2]_0$ &0   
&  $-[\t A_1,A_2]_0$&$[\t v_1,A_2]_0$&$-[\t u_1,A_2]_0$\\                    
& &   &  & $-\frac{1}{2}\langle \t A_1, A_2\rangle_0$  && \\
\hline
$v_2$& 0 &  $-[A_1,v_2]_0$ & 0 & 0  &$0$& 0\\
\hline
$\t A_2$ & 0 &$-[A_1,\t A_2]_0$&  
0 & $\langle\t A_1,\t A_2\rangle_0$& $[\t v_1,\t A_2]_0$&0\\
\hline
$\t v_2$& 0 &   $-[A_1,\t v_2]_0$& 0 & $-[\t A_1,\t v_2]_0 $ &0& 0\\
\hline
$\t u_2$    & $-[A_1,\t u_2]_0 $&   0& 0 & 0  &0& 0\\
\hline
\end{tabular}\\
\vspace{3mm}

In the tables above we keep the same notation as in \rf{mutab} except for 
the fact that 
$u_i v_i, \t u_i, \t v_i$ and $A_i,\t A_i$ are associated with the 
elements from $V^{0}$ and $V^{1}$ correspondingly. 
The bilinear operation $m$ is nonzero only if its both arguments belong to $\mF_1$:
\begin{eqnarray}
m_0(A_1,A_2)=-\langle A_1, A_2\rangle_0. 
\end{eqnarray}
The expression $n_0(a_1,a_2,a_3)$ is nonzero only when all three elements belong to $\mF_1$ or one of the first 
two belongs to $\mF_2$ and the other lie in $\mF_1$:
\begin{eqnarray}
&& n_0(A_1,A_2,A_3)= A_2\langle A_1,A_3\rangle_0-A_1\langle A_2,A_3\rangle_0,\nonumber\\
&& n_0(A_1,\t v,A_2)=n_0(\t v,A_1, A_2)=-\t v \langle A_1, A_2\rangle_0 .
\end{eqnarray}
We believe that the quasiclassical Lian-Zuckerman homotopy BV algebra is actually a $BV_{\infty}$-algebra from  \cite{vallette}, motivated by the results of \cite{gorbounov}, \cite{vallette}. Below we prove an important part of this conjecture, which is needed in the paper.\\

\noindent {\bf Theorem 3.2.} {\it The homotopy associative algebra with the operations $Q,\mu_0, n_0$ 
is an $A_{\infty}$-algebra where the higher homotopies vanish starting from tetralinear one. In other words,  
$Q,\mu_0, n_0$ generate $A_{(3)}$ algebra.}\\

\noindent{\bf Proof.} In order to prove this statement, it is enough to show that the relations:
\begin{eqnarray}
&&(-1)^{n_{a_1}}\mu_0(a_1,n_0(a_2,a_3,a_4))+\mu_0(n_0(a_1,a_2,a_3),a_4)=\\
&&n_0(\mu_0(a_1,a_2),a_3,a_4)-n_0(a_1,\mu_0(a_2,a_3),a_4)+n_0(a_1,a_2,\mu_0(a_3,a_4)),\nonumber
\end{eqnarray}
hold. The last relation is trivial. Let us prove the first one in the most nontrivial case, when $a_i=A_i\in \mF_1$. 

\begin{eqnarray}
&&(-1)^{n_{A_1}}\mu_0(A_1,n_0(A_2,A_3,A_4))+\mu_0(n_0(A_1,A_2,A_3),A_4)=\nonumber\\
&&-\frac{1}{2}\langle A_1,A_3\rangle_0\langle A_2,A_4\rangle_0+\frac{1}{2}\langle A_1,A_2\rangle_0\langle A_3,A_4\rangle_0+\nonumber\\
&&\frac{1}{2}\langle A_2,A_4\rangle_0\langle A_1,A_3\rangle_0-\frac{1}{2}\langle A_1,A_4\rangle_0\langle A_2,A_3\rangle_0=\nonumber\\
&&\frac{1}{2}\langle A_1,A_2\rangle_0\langle A_3,A_4\rangle_0-\frac{1}{2}\langle A_1,A_4\rangle_0\langle A_2,A_3\rangle_0=\nonumber\\
&&n_0(\mu_0(a_1,a_2),a_3,a_4)-n_0(a_1,\mu_0(a_2,a_3),a_4)+n_0(a_1,a_2,\mu_0(a_3,a_4))
\end{eqnarray}
We leave for the reader to establish these relations in all other situations.  
$\hfill \blacksquare$\\

\noindent {\bf 3.3. Remarks on quasiclassical limits and Vertex/Courant algebroids.} In the end of this section we want to draw an attention to an important notion, introduced in \cite{malikov}, called $vertex$ $algebroid$. It basically reflects the relations between the operator products of the elements of conformal dimension $0$ and $1$ and the operator $L_{-1}$ in the vertex algebra (the authors of \cite{malikov} do not require the presence of the Virasoro element). As one can see from the explicit construction of the operations, those relations are included in the relations of 
the LZ algebra \rf{lzrel}. At the same time, for LZ algebra one needs 
extra relations, involving $L_{1}$ operator, which are absent in the definition of vertex algebroid (see e.g. \cite{malikov}, \cite{bressler}, \cite{malikovlag}). 

We also indicate here that the quasiclassical limit of the LZ algebra is not equivalent to the usual quasiclassical limit, which truncates vertex algebra into Poisson algebra, and vertex algebroid into Courant algebroid (see e.g. \cite{bressler}).  In this limit the $L_1$ action vanishes, partly destroying the differential and almost all structures in the homotopy BV algebra, which we obtained, leaving only the subalgebra of the $L_{\infty}$-algebra corresponding to the states of ghost number 0 and 1. This is, in fact, the usual $L_{\infty}$-algebra of the Courant algebroid \cite{weinstein}. 

In the same way the authors of \cite{malikov} build a functor assigning to the vertex algebra the vertex algebroid and therefore $L_{\infty}$-algebra, LZ construction yields a functor which assigns to the 
VOA the homotopy BV algebra. 
Considering VOA with a parameter, we see that we have a functor which maps each such VOA in the quasiclassical limit of LZ algebra, which, as we have seen are highly truncated homotopy BV algebra (see Theorem 3.2). 
As a direct consequence of the Theorem 3.2 we have a functor from the category of VOA with a parameter into the  category of $A_{\infty}$-algebras. 

One might wonder if this construction of the classical limit could be modified to be applied to other TVOA, e.g. the chiral de Rham complex \cite{cdr} in order to obtain some nontrivial homotopy BV structure. However, the space of light modes there is just spanned by differential forms and the operation $\mu$ is given by the wedge product, while 
the analog of $b_0$ acts as a derivation of $\mu$. Therefore the bracket $\{\cdot, \cdot\}$ vanishes and the whole LZ algebra is just the differential graded algebra of differential forms. 
There is, in fact, a nontrivial modification of the LZ algebra in this case, which is worth mentioning here. It basically corresponds to the interchanging of the roles of $b(z)$ and $J(z)$ operators (see Definition 2.1), such that the bracket is constructed by means of $Q$, but not $b_0$. Reducing it to certain subcomplex of $b_0$ one can construct a BV algebra of polyvector fields on a manifold \cite{malikovlag}. This might lead to the chiral version of Barannikov-Kontsevich results \cite{baran}.

\section{Deformation of Lian-Zuckerman homotopy \\ 
commutative algebra.}
{\bf 4.1. Flat background deformation: the case of general TVOA.}
Let $V$ be a TVOA. Let us consider the set of elements $\{f_{i}\}_0$, where
$f_i\in V^0[1]$ are primary elements (i.e. they correspond to the highest weight vectors of the Virasoro algebra)
of conformal dimension 0 and of fermionic degree 1, such that
\begin{eqnarray}
b_0f_i=0,\quad Qf_i=0, \quad \mu(f_i,f_j)=0 \quad \forall i,j.
\end{eqnarray}
An immediate consequence is that $\{f_i,f_j\}_0=0$. 
We will call the operator 
\begin{eqnarray}\label{r}
R^{\eta}=\sum_{i,j}\eta^{ij}\mu(f_i,\{f_j, \cdot\}),
\end{eqnarray}
where $\eta^{ij}$ is some constant matrix, the {\it flat background} deformation of $Q$. 
First of all we show that our definition is consistent, i.e. the following Proposition holds. \\

\noindent {\bf Proposition 4.1.} {\it The operator $R^{\eta}$ obeys the  properties:
\begin{eqnarray}
{(R^{\eta})}^2=0,\quad [Q,R^{\eta}]=0.
\end{eqnarray}}
\noindent{\bf Proof.} The second relation is an immediate consequence of the fact that $Q$ is a derivation of both $\mu(\cdot,\cdot)$ and $\{\cdot,\cdot\}_0$. Let us prove the first one. In order to do that let us write in detail ${(R^{\eta})}^2a$ for some 
$a\in V$. 
\begin{eqnarray}
&&{(R^{\eta})}^2 a=\sum_{i,j,k,l }\eta^{ij}\mu(f_i,\{f_j, \eta^{kl}\mu(f_k,\{f_l, a\})\})=\\
&&\sum_{i,j,k,l}\eta^{ij}\eta^{kl}\mu(f_i\mu(f_k,\{f_j\{f_l, a\}_0\}_0)=\nonumber\\
&&\sum_{i,j,k,l}\eta^{ij}\eta^{kl}(\mu(\mu(f_i,f_k), \{f_j\{f_l, a\}\})+(Qn+nQ)(f_i,f_k, \{f_j\{f_l, a\}\}))\nonumber.
\end{eqnarray}
Since $\{f_j,\{f_l, a\}\}=\{f_l,\{f_j, a\}\}$, ${(R^{\eta})}^2 a=0$. $\hfill \blacksquare$\\

\noindent{\bf Remark: Geometric/Physical Meaning of the Flat Metric Deformation.} 
It can been shown that the operator $R^{\eta}$ has the natural geometric
meaning. Let us consider the operator 
$\phi^{(0)}(z,\b z)=\sum_{i,j}\eta^{ij}f_i(z)f_j(\b z)$, where $\b z$ is the complex conjugated variable for $z$,   
and the operator-valued differential forms $\phi^{(2)}=\sum_{i,j}\eta^{ij}[b_{-1},f_i(z)][b_{-1},f_j(\b z)]$, 
$\phi^{(1)}=d\b z\sum_{i,j}\eta^{ij}[b_{-1},f_i(z)]f_j(\b z) -dz\sum_{i,j}\eta^{ij} f_i(z)[b_{-1},f_j(\b z)]$, such that the following descent equations 
are satisfied:
\begin{eqnarray}
Q\phi^{(2)}=d\phi^{(1)}, \quad Q\phi^{(1)}=d\phi^{(0)}, \quad Q\phi^{(0)}=0.
\end{eqnarray}
Let us consider the operator $\t R^{\eta}$ which acts on the elements of TVOA as follows:
\begin{eqnarray}
\tilde R^{\eta} a=P_0\int_{C_{\epsilon,0}}\phi^{(1)}a,
\end{eqnarray}
where $C_{\epsilon,0}$ is the contour around zero of the radius $\epsilon$ and $P_0$ is the projection on 
$\epsilon^0$-component. Counting the powers of $\epsilon$, one can see that $\t R^{\eta}=R^{\eta}$. 

Physically, it means the following. 
Suppose the TVOA $V$ is described by the action $S_0$. 
Let us consider its perturbation by the operator $\phi^{(2)}$, i.e. the perturbed action is $S_{\phi}=S_0-\int\phi^{(2)}$. Therefore, the deformed $BRST$ current is given by $J_{B,\phi}=J_B(z)dz+\phi^{(1)}$, i.e. 
the deformed charge has the form $Q_{\phi}=\int J_{B,\phi}$, but it must involve some regularization, as it usually happens in the quantum theory. The regularization is given by the projection on $\epsilon$-independent part. After examples are considered, it will be quite clear why we call the perturbation of the form 
$\phi^{(2)}=\sum_{i,j}\eta^{ij}[b_{-1},f_i(z)][b_{-1},f_j(\b z)]$ by "flat background". \\

\noindent {\bf 4.2. Deformation of the quasiclassical LZ algebra of light modes.} 
Let $V$ be a VOA and $V\otimes \Lambda$ is the corresponding BRST complex. Then the following proposition 
holds.\\

\noindent {\bf Proposition 4.2.}{\it The flat background deformation of the BRST operator is determined by Abelian vertex subalgebras of V generated by the primary elements $\{s_i\}$ of conformal weight $1$, such that 
\begin{eqnarray}\label{ab}
s_{i}(z)s_{j}(w)\sim 0, \quad \forall i,j,
\end{eqnarray}
where $s_i(z)=\sum_ns_{n,i}z^{-n-1}$, in other words,  $[s_{n,i}, s_{m,j}]=0$ for any $m,n$}.\\

\noindent {\bf Proof.} Let us consider all possible primary fields of conformal weight $0$ and of ghost number $1$ 
in  $V\otimes \Lambda$. They have the form $cs(z)$ and $\p c u(z)$, where $s(z)$, $u(z)$ are 
primary fields of conformal dimensions 1 and 0 correspondingly. The condition that such field should be annihilated 
by $b_0$ leaves us only the quantum fields of the form $cs(z)$. Therefore, the flat background deformation of the BRST operator is determined by the sets of primary fields $s_i$ of conformal weight 1. Finally, the condition \rf{ab} can be obtained from the fact that $\mu(f_i,f_j)=0$, where $f_i$ is the state corresponding to the quantum field $cs_i(z)$.$\hfill\blacksquare$\\
  
\noindent Therefore, one can introduce the elements $f_i$ from $F_{\mathcal{L}_0}$ corresponding to $s_i(z)$ 
of ghost number $1$ of conformal weight $0$. In the case with parameter, let us define the operator $R_{h}^{\eta}\equiv h^{-1}\sum_{i,j}\eta^{ij}\mu(f_i,\{f_j, \cdot\})$. Then we have a proposition.\\

\noindent {\bf Proposition 4.3.}\\ 
\noindent{\it i) The operator $R_h^{\eta}$ commutes with the BRST operator and ${(R_h^{\eta})}^2=0$. 
It acts  on $F_{\mathcal{L}_0}$ as follows:
\begin{eqnarray}
\xymatrixcolsep{30pt}
\xymatrixrowsep{3pt}
\xymatrix{
& \mathcal{V}^{1}\ar[ddddr]\ar[r]^{\Delta} & \mathcal{V}^{1}\ar[ddddr]^{\frac{1}{2}\hat{L}_1}& \\
&& \hat{L}_{-1} &&\\
& \bigoplus & \bigoplus & \\
&&  -\frac{1}{2}\hat L_{1} &&\\
\mathcal{V}^{0}\ar[uuuur]^{\hat L_{-1}}\ar[r]_{-\Delta} & \mathcal{V}^{0}\ar[uuuur] & \mathcal{V}^{0}\ar[r]_{\Delta}  & \mathcal{V}^{0}
} 
\end{eqnarray}
where $\Delta=h^{-1}\sum_{i,j}\eta^{ij}s_{0,i}s_{0,j}$, $\hat{L}_1\cdot=-h^{-1}\sum_{i,j}\eta^{ij}\langle s_i, s_{0,j}\cdot\rangle$, \\
$\hat{L}_{-1}\cdot=-h^{-1}\sum_{i,j}\eta^{ij}s_is_{0,j}\cdot$.\nonumber\\
\noindent ii) On the complex $\mF$ the operator $R_h^{\eta}$ acts in such a way:
\begin{eqnarray}
R^{\eta}_h:\mF\to \mF[h]
\end{eqnarray}
The quasiclassical version of $R_h^{\eta}$, i.e. the operator 
$R_0^{\eta}=\sum_{i,j}\eta^{ij}\mu_0(f_i,\{f_j, \cdot\}_0)$, where $f_i=cs_i$ and $s_i\in V$, acts on $\mF$ invariantly and commutes with $Q$ on $\mF$.}\\

A natural question is that whether one can deform the quasiclassical LZ algebra in such a way that 
all the relations will be satisfied with the differential $Q^{\eta}=Q+R_{0}^{\eta}$. The answer is only partly positive. 
Namely, the following theorem holds.\\

\noindent{\bf Theorem 4.1.} \\
{\it i)There exist a flat metric deformation of the homotopy associative subalgebra of 
quasiclassical LZ homotopy BV algebra, i.e. there exist $\eta$-deformed multilinear maps $\mu_0^{\eta}$, $m_0^{\eta}$, $n_0^{\eta}$ on $\mF$, which together with $Q^{\eta}$ satisfy the relations of homotopy associative algebra. Moreover, $m_0^{\eta}\equiv m_0$ and $n_0^{\eta}\equiv n_0$.\\ 
\noindent ii) The homotopy associative algebra on $\mF$ with operations $Q^{\eta}, \mu_0^{\eta}, n_0^{\eta}$ is an $A_{\infty}$-algebra, such that all multilinear operations vanish starting from the tetralinear one.
}\\

\noindent{\bf Proof.} 
For simplicity of calculations in this proof we assume Einstein summation convention, i.e. when an index variable appears twice in a single term, once in an upper (superscript) and once in a lower (subscript) position, 
it implies that we are summing over all of its possible values.  

In the beginning, let us show how the bilinear operation is deformed. 
In order to do that let us find an obstacle for $R_0^{\eta}$ to be a derivation of $\mu_0$:
\begin{eqnarray}
&&R_0^{\eta}\mu_0(a_1,a_2)=\mu_0(f_i, \mu_0(\{f_j,a_1\}_0,a_2))\eta^{ij}+\mu_0(f_i, \mu_0(a_1,\{f_j,a_2\}_0))\eta^{ij}=\nonumber\\
&&=\mu_0(\mu_0(f_i, \{f_j,a_1\}_0),a_2)+\mu_0(\mu_0(f_i, a_1),\{f_j,a_2\}_0)\eta^{ij}-\nonumber\\
&&(Qn_0+n_0Q)(f_i,\{f_j,a_1\}_0,a_2)\eta^{ij}-(Qn_0+n_0Q)(f_i,a_1,\{f_j,a_2\}_0)\eta^{ij}=\nonumber\\
&&\mu_0(\mu_0(f_i, \{f_j,a_1\}_0),a_2)+(-1)^{|a_1|}\mu_0(\mu_0(a_1,f_i),\{f_j,a_2\}_0)-\nonumber\\
&&(Qn_0+n_0Q)(f_i,\{f_j,a_1\}_0,a_2)\eta^{ij}-(Qn_0+n_0Q)(f_i,a_1,\{f_j,a_2\}_0)\eta^{ij}+\nonumber\\
&&(Qr_0+r_0Q)(f_i,a_1,\{f_j,a_2\}_0)\eta^{ij}=\nonumber\\
&&\mu_0(\mu_0(f_i, \{f_j,a_1\}_0),a_2)\eta^{ij}+(-1)^{|a_1|}\mu_0(a_1, \mu_0(f_i,[f_j, a_2]))\eta^{ij}+ \nonumber\\
&&(Qr_0+r_0Q)(f_i,a_1,\{f_j,a_2\}_0)\eta^{ij}=\nonumber\\
&&\mu_0(R_0^{\eta} a_1, a_2)+(-1)^{|a_1|}\mu_0(a_1, R_0^{\eta} a_2) -(Q\nu_0^{\eta}+\nu_0^{\eta}Q)(a_1,a_2),
\end{eqnarray}
where $r_0(a_1,a_2)=m_0(f_i,a_1)\{f_j,a_2\}_0\eta^{ij}$ and 
\begin{eqnarray}
\nu_0^{\eta}(a_1,a_2)=n_0(f_i,\{f_j,a_1\}_0,a_2)\eta^{ij}-m_0(f_i,a_1)\{f_j,a_2\}_0\eta^{ij}.
\end{eqnarray}
The explicit expression for $\nu_0^{\eta}$ is given in the following table:  

\begin{center}
$\nu_0^{\eta}(a_1,a_2)$=
\end{center}
\begin{tabular}{|l|c|c|c|c|c|r|}
\hline
\backslashbox{$ a_2$}{$a_1$}&          $u_1$ & $A_1$ & $v_1$ &$\t A_1$ &$\t v_1$& $\t u_1$ \\
\hline
$u_2$ &    0    &$\eta^{ij}\langle s_i,A_1\rangle_0\{f_j,u_2\}_0$ & 0 & 0& 0 & 0\\
\hline
$A_2$     &  0 & $\nu^{\eta}(A_1,A_2)$ &0  & 0 &$-\eta^{ij}\langle s_i,A_2\rangle_0\{f_j,\t v_1\}_0$& 0 \\                    
\hline
$v_2$& 0 &  0 & 0 & 0  & 0 & 0  \\
\hline
$\t A_2$ & 0&  0 &0 & 0&0&0\\
\hline
$\t v_2$& 0 & $-\eta^{ij}\langle s_i,A_1\rangle_0\{f_j,\t v_2\}_0$  & 0 & 0 &0& 0\\
\hline
$\t u_2$    & 0 &   0& 0 & 0  &0& 0\\
\hline
\end{tabular}\\
\vspace{3mm}

\noindent where 
\begin{eqnarray}
&&\nu_0^{\eta}(A_1,A_2)=\\
&&-\langle s_i,A_1\rangle_0\{f_j,A_2\}_0\eta^{ij}-s_i\langle\{f_j,A_1\}_0,A_2\rangle_0\eta^{ij}+\{f_j,A_1\}_0\langle s_i,A_2\rangle_0\eta^{ij}.\nonumber
\end{eqnarray}

We see that $[R_0^{\eta},\mu_0]+[Q,\nu^{\eta}]=0$ and if $[R_0^{\eta}, \nu_0^{\eta}]=0$, 
then $Q^{\eta}$ is a derivation of the bilinear operation $\mu_0^{\eta}=\mu_0+\nu_0^{\eta}$. 
One can show it explicitly. We check the most nontrivial situations, i.e. when $a_1\in \mF_1$, $a_2\in \mF_1$, and also cases when $a_2\in \mF_0$, $a_1\in \mF_1$ and 
$a_1\in \mF_0$, $a_2\in \mF_1$. 

So, let $a_1=u\in\mF_0$ and $a_2=A\in\mF_1$. Then 
$\nu_0^{\eta}(u,A)=0$ and $\nu_0^{\eta}(u,R_0^{\eta}A)=0$. At the same time 
\begin{eqnarray}
&&\nu_0^{\eta}(R_0^{\eta}u,A)=\nu_0^{\eta}(s_i\{f_j,u\}_0,A)\eta^{ij}=\nonumber\\
&&(-s_k\langle s_i\{f_l,\{f_j,u\}_0\}_0,A\rangle_0+s_i\{f_k,\{f_j,u\}_0\}_0\langle s_l,A\rangle_0) \eta^{ij}\eta^{kl}=0.
\end{eqnarray}
Now we check the case $a_2=u\in\mF_0$, $a_1=A\in\mF_1$.
The bilinear operation between these two elements is nontrivial, moreover:
\begin{eqnarray}
R_0^{\eta}\nu_0^{\eta}(A,u)=\eta^{ij}\eta^{kl}s_k\{f_l, \langle s_i, A\rangle_0\{f_j,u\}_0\}_0.
\end{eqnarray}
According to the definition of $\nu_0^{\eta}$, $\nu_0^{\eta}(R_0^{\eta}A, u)=0$. At the same time
\begin{eqnarray}
&&(-1)^{|A|}\nu^{\eta}(A,R_0^{\eta}u)=\langle s_k, A\rangle_0\{ f_l,s_i\{f_j,u\}_0\}_0\eta^{ij}\eta^{kl}+\nonumber\\
&&s_i\langle\{f_j,A\}_0,s_k\{f_l,u\}_0\rangle_0\eta^{ij}\eta^{kl}.
\end{eqnarray}
Therefore, $R_0^{\eta}$ is a derivation in this case also. Let us consider the last possibility, when $a_1=A_1\in\mF_1$ and $a_2=A_2\in \mF_1$.
\begin{eqnarray}
&&R_0^{\eta}\nu_0^{\eta}(A_1,A_2)=\frac{1}{2}\langle s_k,\{f_l,\langle s_i,A_2\rangle_0\{f_j,A_1\}_0\}_0\rangle_0\eta^{ij}\eta^{ij}\eta^{kl}-\nonumber\\
&&\frac{1}{2}\langle s_k,\{f_l,\langle s_i,A_1\rangle_0\{f_j,A_2\}_0\}_0\rangle_0\eta^{ij}\eta^{kl},\nonumber\\
&&\nu_0^{\eta}(R_0^{\eta}A_1,A_2)=\frac{1}{2}\langle s_i,A_2\rangle_0\{f_j,\langle s_k,\{f_l,A_1\}_0\rangle_0\}_0\eta^{ij}\eta^{kl},\nonumber\\
&&-\nu_0^{\eta}(A_1,R_0^{\eta}A_2)=-\frac{1}{2}\langle s_i,A_1\rangle_0\{f_j,\langle s_k,\{f_l,A_2\}_0\}_0\rangle_0\eta^{ij}\eta^{kl}.
\end{eqnarray}
Comparing the terms above, we find that $R_0^{\eta}\nu_0^{\eta}(A_1,A_2)=\nu_0^{\eta}(R_0^{\eta}A_1,A_2)-\nu_0^{\eta}(A_1,R_0^{\eta}A_2)$, i.e. $R_0^{\eta}$ is a derivation.

The next step is to show that $\mu_0^{\eta}$ satisfies the homotopy commutativity relation. 
From the table for $\nu_0^{\eta}$ we see that 
\begin{eqnarray}
&&\nu_0^{\eta}(a_1,a_2)-\nu_0^{\eta}(a_2,a_1)=\nonumber\\
&&R_0^{\eta}\nu_0^{\eta}(a_1,a_2)+\nu_0^{\eta}(R_0^{\eta}a_1,a_2)+
(-1)^{n_{a_1}}\nu_0^{\eta}(a_1,R_0^{\eta}a_2).
\end{eqnarray}
Therefore, $m_0^{\eta}\equiv m_0$. The last statement does mean that $\mu_0^{\eta}$ satisfies the homotopy associativity relation. We will show this in the case when all the arguments belong to $\mF_1$, leaving to check the other combinations to the reader. Let $A_1,A_2,A_3\in \mF_1$, then
\begin{eqnarray} 
&&\mu_0(\nu_0^{\eta}(A_1,A_2),A_3)=\frac{1}{2}(\langle A_3,\{f_j,A_1\}_0\rangle_0 \langle s_i,A_2\rangle_0\eta^{ij}-\nonumber\\
&&\langle f_i,A_1\rangle_0\langle \{f_j,A_2\}_0,A_3\rangle_0\eta^{ij}-\langle s_i,A_3\rangle_0\langle\{f_j,A_1\}_0,A_2\rangle_0\eta^{ij},\nonumber\\
&&\nu_0^{\eta}(\mu_0(A_1,A_2),A_3)=\nu_0^{\eta}(-\frac{1}{2}\langle A_1,A_2\rangle_0,A_3)=\nonumber\\
&&\frac{1}{2}\eta^{ij}\langle f_i,A_3\rangle_0\{f_j,\langle A_1,A_2\rangle_0\}_0,\nonumber\\
&&\nu_0^{\eta}(A_1,\mu_0(A_2,A_3))=\nu_0^{\eta}(A_1,-\frac{1}{2}\langle A_2,A_3\rangle_0)=\nonumber\\
&&\frac{1}{2}\eta^{ij}\{f_j,\langle A_2,A_3\rangle_0\}_0\langle s_i,A_1\rangle_0,\nonumber\\
&&\mu_0(A_1,\nu_0^{\eta}(A_2,A_3))=\frac{1}{2}(\langle A_1,\{f_j,A_2\}_0\rangle_0\langle s_i,A_3\rangle_0\eta^{ij}-\nonumber\\
&&\langle s_i,A_2\rangle_0\langle \{f_j,A_3\}_0,A_1\rangle_0\eta^{ij}-\langle s_i,A_1\rangle_0\langle\{f_j,A_2\}_0,A_3\rangle_0\eta^{ij}).
\end{eqnarray}
Therefore, 
\begin{eqnarray} 
&&\mu_0(\nu^{\eta}(A_1,A_2),A_3)+\nu_0^{\eta}(\mu_0(A_1,A_2),A_3)-\nonumber\\
&&\nu_0^{\eta}(A_1,\mu_0(A_2,A_3))-\mu_0(A_1,\nu^{\eta}(A_2,A_3))=\nonumber\\
&&\frac{1}{2}\eta^{ij}\{f_j,\langle A_3,A_1\rangle_0\}_0\langle s_i,A_2\rangle_0-\frac{1}{2}\eta^{ij}\{f_j,\langle A_2,A_3\rangle_0\}_0\langle s_i,A_1\rangle_0.
\end{eqnarray}
On the other hand, 
\begin{eqnarray} 
&&R_0^{\eta}n_0(A_1,A_2,A_3)+n_0(R_0^{\eta}A_1,A_2,A_3)- n_0(A_1,R_0^{\eta}A_2,A_3)+\nonumber\\
&&n(A_1,A_2,R_0^{\eta}A_3)=\nonumber\\
&&\frac{1}{2}\eta^{ij}\{f_j,\langle A_3,A_1\rangle_0 \langle s_i,A_2\rangle_0\}_0-\frac{1}{2}\eta^{ij}\{f_j,\langle A_2,A_3\rangle_0\langle s_i,A_1\rangle_0\}_0+\nonumber\\
&&\frac{1}{2}\langle A_2,A_3\rangle_0\langle s_i,\{f_j,A_1\}_0\rangle_0 \eta^{ij}-\frac{1}{2}\langle A_1,A_3\rangle_0\langle s_i,\{f_j,A_2\}_0\rangle_0 \eta^{ij}=\nonumber\\
&&\mu_0(\nu^{\eta}(A_1,A_2),A_3)+\nu_0^{\eta}(\mu_0(A_1,A_2),A_3)-\nonumber\\
&&\nu_0^{\eta}(A_1,\mu_0(A_2,A_3))-\mu_0(A_1,\nu_0^{\eta}(A_2,A_3)).
\end{eqnarray}
Hence, $\mu_0^{\eta}$ is associative up to homotopy provided by $n$.  Hence, we proved $i)$. In order to prove $ii)$ we just note that $\nu^{\eta}$ doesn't contribute to the higher associativity relation involving 
$\mu_0^{\eta}, n_0^{\eta}$. Then $ii)$ follows from the proof of the similar statement for $Q,\mu_0, n_0$.
$\hfill \blacksquare$\\

Since we have an $A_{\infty}$-algebra it is natural to consider the associated generalized Maurer-Cartan equation. However, due to the properties of the operations $\mu_0, n_0$, the resulting equation coincides with the linear one. In order to get around this, we multiply the  $A_{\infty}$-algebra with some noncommutative associative algebra $S$. The resulting object is also from the category of $A_{\infty}$-algebras, but there will be no homotopy commutativity in general. Let us choose $S=U(\mathfrak{g})$, the universal enveloping algebra of some Lie algebra $\mg$. In this case it is possible (as we will see on several examples below) to find the relation with gauge theory.\\

\noindent{\bf Definition 4.1.}{\it Consider the $A_{\infty}$-algebra on $\mF\otimes U(\mg)$ generated by $Q^{\eta}, \mu_0^{\eta}, n_0^{\eta}$. Consider the Maurer-Cartan element $\Psi\in\mF_1\otimes \mg$. We will refer to the Maurer-Cartan equation
\begin{eqnarray}
Q^{\eta}\Psi+\mu_0^{\eta}(\Psi,\Psi)+n_0^{\eta}(\Psi,\Psi,\Psi)=0
\end{eqnarray}
as the Yang-Mills equation associated to the VOA with a formal parameter $\mathcal{V}$ and Lie algebra $\mg$.
We will refer to the infinitesimal symmetries of the generalized Maurer-Cartan equations, which are 
\begin{eqnarray}
\Psi\to \Psi+\epsilon(Q^{\eta}u+\mu_0^{\eta}(\Psi,u)-\mu_0^{\eta}(u,\Psi)),
\end{eqnarray}
where $u\in\mF_0\otimes\mg$, as gauge symmetries. }\\

\noindent In the next section we will show that this equation and its gauge symmetries are actually 
equivalent to the system of Yang-Mills equations with matter fields and their gauge symmetries for certain vertex algebras.  

In the end of this section we write the explicit expression for the Yang-Mills equation associated with VOA $V$ and the Lie algebra $\mg$, since it is needed in the following. 

The Maurer-Cartan element for $\mF\otimes U(\mg_)$ has the form $\Psi=A+v$, where $A,v\in\mF_1\otimes \mg$ are the elements corresponding to the states of conformal weights $1$, $0$ correspondingly. The equation for $v$ is as 
follows:
\begin{eqnarray}
v=\frac{1}{2}L_1A+\frac{1}{2}\eta^{ij}\langle s_i,\{s_j,A\}_0\rangle_0+\frac{1}{2}\langle A,A\rangle_0.
\end{eqnarray}
Then the equation for $A$ is:
\begin{eqnarray}\label{ymeq}
&&2\Delta A+L_{-1}L_1A+\sum_{i,j}\eta^{ij}L_{-1}\langle s_i,\{s_j,A\}_0\rangle_0+\nonumber\\
&&\sum_{k,l}\eta^{kl}s_k\{s_l,L_1A+\sum_{i,j}\eta^{ij}\langle s_i,\{s_j,A\}_0\rangle_0\}_0+\nonumber\\
&&A(L_1A)+A\sum_{i,j}\eta^{ij}\langle s_i,\{s_j,A\}_0\rangle_0-(L_1A)A+\nonumber\\
&&\sum_{i,j}\eta^{ij}\langle s_i,\{s_j,A\}_0\rangle_0 A+2[A,A]-\langle A^{ad},A^{ad} \rangle_0 A=0,
\end{eqnarray}
where $A^{ad}\in \mF_1\otimes End(\mg)$ stands 
for the element $A\in \mF_1\otimes \mg$ with the Lie algebra elements are considered in the adjoint representation. 
The gauge symmetries of this equation can be written as follows:
\begin{eqnarray}\label{ymsym}
A\to A+\epsilon(L_{-1}u+\eta^{ij}s_i\{s_j,u\}_0+Au-uA),
\end{eqnarray}
where $u\in\mF_0\otimes \mg$. 

\section{ Beta-gamma systems, Courant Algebroid and Yang-Mills Theory}
\noindent {\bf 5.1. A toy model: Heisenberg VOA.} Let us start with the simplest nontrivial vertex algebra with a formal parameter, i.e. we consider the vertex algebra $V(a,g)[h, h^{-1}]$, generated by quantum fields $a^i(z)$ $(i=1,...,D)$ 
with operator products   
\begin{eqnarray}
a^i(z)a^j(w)\sim \frac{hg^{ij}}{z-w},
\end{eqnarray}
where $g^{ij}$ is a symmetric matrix. 
In the case of this VOA the elements of zero conformal weight are just constants and the fields of conformal 
dimension 1 have the form $A_ia^i(z)$, where $A_i$ are constant. This example is very special 
because of the following two facts presented in Proposition 5.1.\\ 

\noindent {\bf Proposition 5.1.} {\it i) The quasiclassical limit of the LZ algebra of light modes on $V(a,g)[h, h^{-1}]$ is isomorphic to the 
LZ algebra of light modes on $V(a,g)$.\\
ii) Any flat background deformation of BRST operator is trivial for $V(a,g)$. 
}\\

\noindent {\bf Proof.} To prove i) it is enough to note that the terms corresponding to higher powers in $h$ usually correspond to the multiple "contractions" of the quantum fields. In the case of $V(a,g)$ this does not happen, since all the elements from the light mode complex are at most linear in $a^i(z)$. 

Expressing $a^i(z)=\sum_na^i_nz^{-n-1}$, we find that $a^i_0$ annihilates each element in $V(a,g)$ for every $i$. Therefore, ii) is also proven. $\hfill \blacksquare$\\

\noindent The next Proposition is about the Yang-Mills equation on $V(a,g)$.\\
 
\noindent {\bf Proposition 5.2.} {\it The Yang-Mills equation for $V(a,g)$ with Lie-algebra $\mg$ 
is equivalent to the system of equations 
\begin{eqnarray}
\sum_{i,j}g^{ij}[A_i,[A_j,A_k]]=0
\end{eqnarray}
for certain $A_i\in \mg$ $(i=1,..., D)$}

One can see that these equations coinside with Yang-Mills equations with a flat metric $g^{ij}$ and constant gauge fields $A_i$.

However, this example was a toy model for us: it was too degenerate. In order to have less trivial example, one should "enrich" the space of fields of conformal weight zero. \\

\noindent {\bf 5.2. BV extension of the Courant algebroid.} 
Let us consider a family of $\beta$-$\gamma$ systems generated by quantum fields 
$p_i(z)$, $X^i(z)$, $(i=1,...,D)$, where $p_i(z)=\sum_n p_{i,n}z^{-n-1}$, $X^i(z)=\sum_n X^i_nz^{-n-1}$ are quantum fields of conformal dimensions $1$ and $0$ correspondingly, and the operator products are
\begin{eqnarray}
X^{i}(z)p_{j}(w)\sim\frac{h\delta_{i,j}}{z-w}, \quad X^{i}(z)X^{j}(w)\sim 0, \quad p_{i}(z)p_{j}(w)\sim 0
\end{eqnarray}
and the Virasoro element is given by the formula
\begin{eqnarray}
L(z)=-\frac{1}{h}\sum_ip_i\p X^i.
\end{eqnarray}
The space of the VOA is given by $F_{p,X}[h^{-1},h]$, $F_{p,X}\equiv\otimes^D_{i=1}F_{p_i,X^i}$ (see section 3.1). 
For definiteness, from now on let us assume that the space of zero conformal weight of $F_{p,X}$ is given by the formal powers series in $X_0^i$. The operators from the VOA $F_{p,X}$ of conformal dimensions $0$ and $1$ have the form 
\begin{eqnarray}
&&u(z)=u(X)(z),\quad  \mA(z)=\sum_i:A^{i}(X)(z)p_{i}(z):, \nonumber\\
&&\mB(z)=\sum_jB_{j}(X)(z)\p X^{j}(z),
\end{eqnarray}
where $u(X), A^{i}(X), B_{j}(X)$ are considered as power series in $X^i$. Let $M$ be the formal scheme 
${\rm{Spf}}(\mathbb{C}[[X_0^1,...,X_0^D]])$. Then the states $u$, $\mA$, $\mB$ can be identified with sections of $\mathcal{O}_M$, $TM$, $T^*M$ correspondingly. 

Now let us consider the semi-infinite complex associated with $F_{p,X}$. The BRST operator acts as follows:

\begin{eqnarray}\label{qbg}
\xymatrixcolsep{40pt}
\xymatrixrowsep{5pt}
\xymatrix{
& TM\ar[ddr]^{\frac{1}{2}\rm{div}} & TM\ar[ddr]^{-\frac{1}{2} \rm{div}} & \\
& \bigoplus & \bigoplus & \\
\mathcal{O}_M\ar[ddr]_{\rm{d}} & \mathcal{O}_M\ar[r]^{\rm{id}}\ar[ddr]_{\rm{d}} & \mathcal{O}_M &\mathcal{O}_M\\
& \bigoplus & \bigoplus &\\
& T^*M & T^*M &
}
\end{eqnarray}
where ${\rm div} \mA=\sum_i\p_{i}A^{i}$. The action of the BV operator $b_0$ is given by the diagram below.

\begin{eqnarray}
\xymatrixcolsep{40pt}
\xymatrixrowsep{5pt}
\xymatrix{
& TM & TM\ar[l]_{-\rm{id}} & \\
& \bigoplus & \bigoplus & \\
\mathcal{O}_M & \mathcal{O}_M\ar[l]_{\rm{id}} & \mathcal{O}_M &\mathcal{O}_M\ar[l]_{-\rm{id}}\\
& \bigoplus & \bigoplus &\\
& T^*M & T^*M\ar[l]_{-\rm{id}} &
}
\end{eqnarray}
It is useful to write down the explicit values for the 
operation $\mu$ in the LZ algebra on the complex $\mF$:

\begin{center}
$\mu_0(a_1,a_2)$=
\end{center}
\begin{tabular}{|l|c|c|c|c|c|r|}
\hline
 \backslashbox{$a_2$}{$a_1$}&          $u_1$ & $\mathcal{X}_1$ & $v_1$ &$\t {\mathcal{X}_1}$ &$\t v_1$& $\t u_1$ \\
\hline
$u_2$ &                  $u_1u_2$    &$\mathcal{X}_1 u_2$ &$v_1u_2$&
$\t {\mathcal{X}_1} u_2$&$\t v_1 u_2$ &$\t u_1 u_2$\\
&  & $-L_{\mathcal{X}_1}u_2$  &  &   && \\
\hline
$\mathcal{X}_2$     &  $u_1\mathcal{X}_2$ & $(\mathcal{X}_1,\mathcal{X}_2)_D+$ &$-v_1\mathcal{X}_2$  
&  $-\frac{1}{2}\langle \t {\mathcal{X}_1},\mathcal{X}_2\rangle$&$L_{\mathcal{X}_2}\t v_1$&0\\                    
& & $\frac{1}{2}\langle \mathcal{X}_1, \mathcal{X}_2\rangle$  &  &  && \\
\hline
$v_2$& $u_1\t u_2$ &  $v_2\mathcal{X}_1$ & 0 & 0  &$-\t v_1v_2$& 0\\
\hline
$\t {\mathcal{X}_2}$ & $u_1\t {\mathcal{X}_2}$ &$-\frac{1}{2}\langle \mathcal{X}_1,\t {\mathcal{X}_2}\rangle$&  
0 & 0& 0&0\\
\hline
$\t v_2$& $u_1\t u_2$ &   $L_{\mathcal{X}_1}\t v_2$& $-v_1\t v_2$ & 0  &0& 0\\
\hline
$\t u_2$    & $u_1\t u_2$ &   0& 0 & 0  &0& 0\\
\hline
\end{tabular}\\
\vspace{3mm}


\noindent where $\mathcal{X}_i,\t {\mathcal{X}_i}$ stand for pairs $(\mA_i,\mB_i), (\t \mA_i,\t \mB_i)$, such that $\mA_i,\t \mA_i\in TM$ and 
$\mB_i,\t \mB_i\in T^*M$.  
The expression $(\mathcal{X}_1,\mathcal{X}_2)_D$ stands for the Dorfman bracket (see e.g. \cite{courant}, \cite{roytenberg}) defined as 
\begin{eqnarray}
((\mA_1,\mB_1), (\mA_2,\mB_2))_D=(L_{\mA_1}\mA_2, L_{\mA_1}\mB_2-i_{\mA_2}d\mB_1).
\end{eqnarray}
The pairing $\langle \cdot,\cdot \rangle$ between two elements $(\mA_1,\mB_1)$ and $(\mA_2,\mB_2)$ is defined 
as follows:$\langle (\mA_1,\mB_1), (\mA_2,\mB_2)\rangle=\sum_{\mu}(A_1^{\mu}B_{2,\mu}+B_{1,\mu}A_{2}^{\mu})$, i.e. it  
is a natural pairing between the elements of $TM$ and $T^*M$. The operation $L_\mathcal{X} u$ for $\mathcal{X}=(\mA,\mB)$  
denotes the Lie derivative with respect to the vector field $\mA$. 

This suggests that there is a certain relation between the Courant algebroid \cite{courant} on 
$TM\oplus T^*M$ and the quasiclassical LZ algebra for $\beta$-$\gamma$ systems. Namely the following Proposition holds.\\

\noindent {\bf Proposition 5.3.} {\it The homotopy BV algebra generated by 
$\mu_0$ and $\{,\}_0$ on $\mF$ contains the Courant algebroid structure on $TM\oplus T^*M$, i.e. 
 $\{,\}_0$ being restricted to $TM\oplus T^*M$ coincides with the Dorfman bracket, $\mu_0: \mF_0\otimes\mF_1\to \mF_1$ is a multiplication of function on an element of $TM\oplus T^*M$, and the pairing between two elements of 
 $\mF_1$ is given by the operator product coefficient $\langle\cdot,\cdot \rangle$.}\\

\noindent In fact, this relation between such "short" homotopy BV algebras and Courant brackets 
can be extended to the large class of Courant algebroids  (we will discuss this and related questions in \cite{bvcourant}). 

In the following we will call this homotopy BV algebra the $BV$ $double$ of Courant algebroid.\\

\noindent {\bf 5.3. Deformation of the homotopy commutative algebra and Yang-Mills equations.} In the previous subsection we showed that the quasiclassical Lian-Zuckerman algebra associated with a family of $\beta$-$\gamma$ systems gives a homotopy $BV$ algebra extending the Courant algebroid. In this section we consider the flat background deformation of this algebra according to considerations of the previous section. 
First of all we will pick the Abelian subalgebra of the operators of conformal dimension 1 in the beta-gamma VOA. We  will consider the VOA subalgebra generated by $p_i(z)$, $(i=1,...,D)$. Moreover, in this subsection 
we assume that the deformation matrix $\eta^{ij}$ is symmetric. Then the deformation of the BRST operator is given in the following proposition.\\

\noindent {\bf Proposition 5.4.} {\it Let $f_i(z)=cp_i(z)$. Then $R_0^{\eta}=\sum_{\i,j}\eta^{ij}\mu_0(f_i,\{f_j,\cdot\}_0)$ acts on $\mF$ as follows:
\begin{eqnarray}
\xymatrixcolsep{40pt}
\xymatrixrowsep{5pt}
\xymatrix{
& TM\ar[r]^\Delta & TM & \\
& \bigoplus & \bigoplus & \\
\mathcal{O}_M\ar[uur]^{\hat{\rm{d}}}\ar[r]^{-\Delta} & \mathcal{O}_M\ar[uur]^{\hat{\rm{d}}} & \mathcal{O}_M\ar[r]^\Delta & \mathcal{O}_M\\
& \bigoplus & \bigoplus &\\
& T^*M\ar[uur]^{\frac{1}{2}\widehat{\rm{div}}} \ar[r]^\Delta & T^*M\ar[uur]^{-\frac{1}{2}\widehat{\rm{div}}}&
}
\end{eqnarray}
where $\Delta=\sum_{i,j}\eta^{ij}\p_i\p_j$, $\widehat{\rm{div}}\mB=\sum_{i,j}\eta^{ij}\p_iB_j$ (here $\mB\in T^*M$) and $(\hat{\rm{d}}u)^j=\sum_{i}\eta^{ij}\p_i u$ (here $u\in \mathcal{O}_M$).} \\

\noindent Let us also assume that the matrix $\eta^{\{ij\}}$ is invertible, such that $\eta_{\{ij\}}$ is the inverse matrix. 
This yields the following proposition.\\

\noindent {\bf Proposition 5.5.}{\it The complex $(\mF, Q^{\eta})$ is isomorphic to the following 
complex $(\mG, {Q'}^{\eta})$ which decomposes into a direct sum of three subcomplexes: 
\begin{eqnarray}\label{omegacomp}
\xymatrixcolsep{30pt}
\xymatrixrowsep{-5pt}
\xymatrix{
0\ar[r]&\Omega^0(M) \ar[r]^{\ud} & \Omega^1(M) \ar[r]^{*\ud*\ud}  &\Omega^1(M) \ar[r]^{*\ud*} & \Omega^0(M)\ar[r]&0\\
 &\quad &     & \quad    &&\\
 && \bigoplus & \bigoplus     &&\\
 &&           & \quad    &&\\
 &0\ar[r]& \Omega^1(M) \ar[r]^{\Delta}   & \Omega^1(M)\ar[r]&0\\
&& \bigoplus & \bigoplus     &&\\
&0\ar[r]& \Omega^0(M) \ar[r]^{id}   & \Omega^0(M)\ar[r]&0
}
\end{eqnarray}
where $\Omega^0(M)\equiv \mathcal{O}_M$ and $\Omega^1(M)\equiv T^*M$ and the Hodge star operator $*$ is constructed 
via the metric corresponding to the invertible and symmetric matrix $\eta_{\{ij\}}$.}\\

\noindent{\bf Proof.} The embedding can be constructed in the following way. Let us denote the three 
subcomplexes above as $(\mathfrak{G}_i^{\cdot},{Q'}^{\eta})$ $(i=1,2,3)$,   
We construct the following emebeddings:
\begin{eqnarray}
&&\mathfrak{G}_1^0\xrightarrow{id} \mF_0, \quad \mathfrak{G}_1^3\xrightarrow{id} \mF_3,\nonumber\\
&&\mathfrak{G}_1^1\xrightarrow{f_1}\mF_1, \quad  \mathfrak{G}_1^1\xrightarrow{g_1}\mF_1,\nonumber\\
&&\mathfrak{G}_2^1\xrightarrow{f_2}\mF_1, \quad  \mathfrak{G}_2^1\xrightarrow{g_2}\mF_1,\nonumber\\
&&\mathfrak{G}_3^1\xrightarrow{f_3}\mF_1,  \quad \mathfrak{G}_3^1\xrightarrow{g_3}\mF_1, 
\end{eqnarray}
such that $f_i, g_i$ are explicitly given by:
\begin{eqnarray}
&&f_1(\mB)=\mB+\mB^*-\hat{{\rm div}}\mB,\quad  g_1(\t \mB)=\t \mB+\t \mB^*,\nonumber\\
&&f_2(\mB)=\mB-\mB^*, \quad g_2(\t \mB)=\t \mB-\t \mB^*,\nonumber\\
&&f_3\equiv id, \quad g_3(\t v)=\t v- (d\t v+\hat d\t v), 
\end{eqnarray}
where $\mB^*\in TM$ such that ${B^*}^i=\eta^{ij}B_j$. Combining $f_i, g_i$ and other maps into the map of complexes, one can see that this is an isomorphism.
$\hfill\blacksquare $ \\

It can be observed that the cohomology of the complex $(\mF,Q^{\eta})$ in degree 1 is equivalent to the solutions of Maxwell equations and scalar field equations
\begin{eqnarray}
*\ud*\ud\mA=0, \quad \Delta\Phi=0,
\end{eqnarray} 
modulo gauge transformations $\mA\to \mA+\ud u$. Now we recall that we have the structure of the homotopy associative and homotopy commutative algebra on $(\mF,Q^{\eta})$, which is in fact the $A_{\infty}$-algebra. 
Let us consider the tensor product of $(\mF,Q^{\eta})$ with $U(\mathfrak{g})$, where $\mathfrak{g}$ is some Lie algebra and find out what are the Yang-Mills equations for the $\beta$-$\gamma$ VOA. 
The answer is given in the proposition below.\\

\noindent {\bf Proposition 5.6.} {\it The Yang-Mills equation for the VOA $F_{p,X}$ and the Lie algebra $\mg$ 
is equivalent to the following system of equations
\begin{eqnarray}\label{aphi}
&&\sum_{i,j}\eta^{ij}[\nabla_{i},[\nabla_{j},\nabla_{k}]]=\sum_{i,j}\eta^{ij}[[\nabla_{k},\Phi_{i}],\Phi_{j}],\nonumber\\
&&\sum_{i,j}\eta^{ij}[\nabla_{i},[\nabla_{j},\Phi_{k}]]=\sum_{i,j}\eta^{ij}[\Phi_{i},[\Phi_{j},\Phi_{k}]],
\end{eqnarray}
where $\Phi_{i}={B}_{i}- \sum_jA^{j}\eta_{ij}$, $\mathcal{A}_{i}=B_{i}+\sum_jA^{j}\eta_{ij}$ and $\nabla_{i}=\p_{i}+\mathcal{A}_{i}$.  Here $B_i$ are the components of $\mB\in T^*M\otimes \mg$ and $A_i$ are the components of $\mA\in TM\otimes \mg$.\\
The gauge symmetries correspond to the following transformation of fields:
\begin{eqnarray}\label{aphisym}
\mathcal{A}_{i}\to \mathcal{A}_{i}+\epsilon(\p_i u+[\mathcal{A}_{i},u]), \quad 
\Phi_{i}\to \Phi_{i}+\epsilon[\Phi_{i},u].
\end{eqnarray}
}
\noindent {\bf Proof.} To prove this statement we just need to substitute the generic Maurer-Cartan element, which is the sum of $(\mA,\mB, v)\in (TM\oplus T^*M\oplus\mathcal{O}_M)\otimes\mg$. Substituting it in \rf{ymeq}, we obtain the equations $\rf{aphi}$. The equation for $v$ is as follows:
\begin{eqnarray}
v=-\frac{1}{2}(\sum_i\p_iA^i+\sum_{ij}\eta^{ij}\p_iB_{j})-\frac{1}{2}\sum_i(A_iB^i+B^iA_i).
\end{eqnarray}
Substituting appropriate elements in \rf{ymsym}, we obtain gauge symmetries \rf{aphisym}. $\hfill \blacksquare$ \\

\noindent 
This statement is very close to the particular results obtained in relation to the "original" logarithmic open string vertex algebra \cite{bvym}, \cite{cftym}. The equations \rf{aphi} are the Yang-Mills equations in the presence 
of $D$ scalar fields $\Phi_i$ (see also \cite{azbg}).  
\\
%

\noindent {\bf 5.4. Smooth manifold case.}
All the considerations we had above were applied only to the case of flat metric and a standard volume form on $D$-dimensional (pseudo-)Euclidean space. Here we will give some statements about the case of a general smooth manifold $M$. First of all we generalize the  BV double of the Courant algebroid to the case of a $D$-dimensional smooth manifold $M$ with a volume form $\Omega$, such that in local coordinates $\Omega=e^{\phi(X)}dX^1\wedge....\wedge dX^D$. 
We have the following proposition.\\

\noindent {\bf Proposition 5.7.} {\it Let us consider the quasiclassical LZ algebra for the VOA 
$F_{p,X}$, such that the Virasoro element is given by 
\begin{eqnarray}
L^{\phi}(z)=-\frac{1}{h}\sum_ip_i\p X^i+\p^2\phi(X).
\end{eqnarray}
Then there exist a BV algebra defined on the sections of certain bundles of the manifold $M$, such that in local coordinates it is given by this quasiclassical LZ algebra.}\\

\noindent {\bf Proof.}  The shift in  $\p^2\phi(X)$ of the Virasoro element changes the action of the differential in \rf{qbg} in such a way that the ${\rm div}$ operation is changed by ${\rm div}\mA=\sum_i(\p_iA^i+\p_i\phi A^i)$. This is a local coordinate expression for the operator invariant under the coordinate change. As one can see, the other operations are already written in the covariant form, therefore, we have a BV algebra defined globally on the sections of appropriate bundles on $M$. $\hfill \blacksquare$\\

\noindent The next question is that whether the deformed homotopy commutative algebra can be generalized to the case of the manifold $M$ with the metric $g=g_{ij}dX^idX^j$. The answer is positive and the following Proposition holds.\\

\noindent {\bf Proposition 5.8.} {\it There exists an $A_{\infty}$-algebra on the complex \rf{omegacomp} for the 
smooth Riemannian manifold $(M,g)$, such that in the case of $M=\mathbb{R}^D$ and $g^{ij}=\eta^{ij}$ it coincides with the $A_{\infty}$-algebra, generated by $Q^{\eta},\mu^{\eta},n^{\eta}$.}\\

\noindent {\bf Proof.} Similar $A_{\infty}$-algebra for the general manifold $M$  
on the sum of upper and lower complexes from \rf{omegacomp} was considered in \cite{cftym}. In the same way one can construct it on the full complex $\mF$. We leave the technical details to the reader.  $\hfill \blacksquare$\\

\noindent For a compact manifold $M$, one can write the action for the Yang-Mills theory interacting with 1-forms, corresponding to this $A_{\infty}$-algebra:
\begin{eqnarray}
&&S=\int_M (-\frac{1}{4}\mathbf{F}\wedge *\mathbf{F}+\frac{1}{4}d_{\mathcal{A}}\Phi\wedge *d_{\mathcal{A}}\Phi+  \frac{1}{4}d_{\mathcal{A}}*\Phi\wedge *d_{\mathcal{A}}*\Phi +\nonumber\\
&&\frac{1}{2}\mathbf{F}\wedge *(\Phi\wedge\Phi)-\frac{1}{4}(\Phi\wedge\Phi)\wedge *(\Phi\wedge\Phi)),
\end{eqnarray}
where $\mathbf{F}=d\mathcal{A}+\mathcal{A}\wedge\mathcal{A}$ is a curvature for $\mathcal{A}$ 
and $d_{\mathcal{A}}=d+\mathcal{A}$ is a covariant derivative. 

It is not clear, however, how to derive the Yang-Mills equations on the smooth manifold $M$ with some metric $g_{ij}$
from a point of view of VOA, like we did before in the flat case (see some suggestions in the last section).
However, the following statement is still true. \\

\noindent {\bf Proposition 5.9.} {\it Let us consider the $A_{\infty}$-algebra on the manifold $M$ with the metric $g$ 
from the Proposition 5.8. Let us introduce a formal parameter into $g$, such that $g^{ij}\to g^{ij}(t)=t g^{ij}$. Then 
taking the limit $t\to 0$ we recover the $A_{\infty}$-subalgebra of the BV Courant algebroid on M with a volume form 
$\Omega=\sqrt{g(X)}dX^1\wedge...\wedge dX^N$, where $g(X)=det(g_{\{ij\}}(X))$. }\\

\noindent{\bf Proof.}  In the flat case this result is obvious. For nonconstant $g^{ij}$ we need to watch 
that terms containing $\sqrt{g(t)}$ will not blow up in the $t\to 0$ limit. This never happens, since they always enter the expressions in the form $\p_i\log(g(t))\equiv\p_i\log(g)$. $\hfill \blacksquare$

\section{Some remarks, open problems and conjectures}
\noindent {\bf 6.1. Physical interpretation: beta-functions and (deformed) LZ algebras.}
In this paper we have shown that the correspondence between VOA and $A_{\infty}$-algebras we constructed using quasiclassical limit of the Lian-Zuckerman homotopy algebra of light modes 
allows us to write down an analogue of the Yang-Mills equation for the general VOA with a formal parameter. It is known that Yang-Mills equations show up as a 1-loop beta-function for the open string theory. Therefore, one can 
give the pure algebraic meaning to the 1-loop beta function in the generic case. 
In such a way the Maurer-Cartan element is associated with a perturbation term. One can think of that as follows.  
Let the VOA correspond to some CFT with the action $S_0$ on the half-plane $H^+$. Then the flat background term may be interpreted as a perturbation of an action of the form $\int_{H^+}\frac{1}{h}\sum_{i,j}\eta^{ij}s_i(z)s_j(\b z)$ as we already have 
seen. The perturbation corresponding to the Maurer-Cartan element $\Psi$ has the form $\int_{\p H^+}b_{-1}\Psi$. 
Therefore, the complete action of the theory (sigma-model) for which we write the conformal invariance condition has the form (for $\beta$-$\gamma$ example such action was considered in \cite{azbg}):
\begin{eqnarray}
S=S_0+\int_{H^+}\frac{1}{h}dz\wedge d\b z (\sum_{i,j}\eta^{ij}s_i(z)s_j(\b z))+\int_{\p H^+}dz (b_{-1}\Psi).
\end{eqnarray}
In the case of non-Abelian Lie algebra (when $\Psi\in \mF\otimes\mg$) one cannot just add the last term to the 
action. One has to insert a trace of the path ordered exponential of $b_{-1}\Psi$ in the partition function. 

One of the open mathematical questions, which is motivated by this paper, is to 
prove partly the Lian-Zuckerman conjecture, i.e. reconstruct completely the LZ algebra of light modes, and to compare the corresponding Maurer-Cartan equations 
with the known expressions for the beta-functions. 
In particular, we partly know the explicit form of the Maurer-Cartan equations in certain nontrivial situations, i.e. in the standard open string case, the low derivative terms should correspond to the equation in the Born-Infeld theory.
  
One of the ways to construct this homotopy algebra is to treat it as a deformation of the quasiclassical one. However, there can be several such deformations. If so, it is also interesting to know what does it mean from the physical point of view, i.e. from the point of view of perturbation theory. \\

\noindent {\bf 6.2. $B$-field and the deformations of the BV homotopy algebra.} In the example, studied 
in subsection 5.2, corresponding to the vertex algebra, generated by 
$\beta$-$\gamma$ systems for simplicity we considered the deformation via the symmetric matrix $\eta^{ij}$. In the case of general $\eta$, the antisymmetric part is related to the so-called Kalb-Ramond field which is a necessary ingredient in string theory. It is interesting to write the Yang-Mills equations and the corresponding action in this instance.  

 Another interesting structure arises when the matrix $\eta^{ij}$ becomes purely antisymmetric. We state here a proposition to which we will return in the subsequent publications.\\
 
\noindent {\bf Proposition 6.1.}  {\it Let $\eta^{ij}$ be antisymmetric. Then there exists a homotopy BV algebra on $\mF$ such that 
it is generated by $Q^{\eta},\mu_0^{\eta},  \t b_0=h^{-1}b_0, \{\cdot ,\cdot \}^{\eta}_0$, 
such that $\{a_1 ,a_2 \}^{\eta}_0=\t b_0\mu_0^{\eta}(a_1,a_2)-
\mu_0^{\eta}(\t b_0a_1,a_2)-(-1)^{n_{a_1}}\mu_0^{\eta}(a_1,\t b_0a_2)$. }\\

We note, that if $\eta^{ij}$ is not antisymmetric, such BV algebra does not exist, the obstacle corresponds to the generalization of the Laplace operator $\Delta\cdot=\sum_{i,j}\eta^{ij}\{s_i,\{s_j,\cdot\}$. 

 As the simplest example, let us consider $\beta$-$\gamma$ system. In this case the deformation of the corresponding homotopy BV algebra via antisymmetric tensor (related to Abelian vertex subalgebra generated by $p_{\{i\}}$) leads to a certain deformation of the BV double of Courant algebroid via the antisymmetric bivector field. One can check that the corresponding deformation of the Courant algebroid coincides with the one considered in \cite{roytenberg} (see also \cite{halmagyi} for the physical 
insight). It appears that a necessary condition on 
bivector field $\eta$ for this deformation to hold on general smooth manifold $M$ 
was that $[\eta,\eta]_S=0$, where $[\cdot,\cdot]$ is a Schouten bracket. It is known, that the Courant algebroid 
admits also another deformation (or a twist), via the 3-form (see e.g. \cite{roytenberg}) which was introduced by Severa. It is interesting, how to incorporate this deformation in the VOA formalism we studied in this paper. One of the ways to do that is to consider instead of $\eta$-deformation based on the Abelian vertex algebra generated by $p_i$, 
another deformation based on the abelian vertex algebra generated by vertex 
operators of the type $\omega_k(X)\p X^k$.
\\

\noindent {\bf 6.3. Yang-Mills equations for $\beta$-$\gamma$ systems and T-duality.} Looking at the Yang-Mills equations in the $\beta$-$\gamma$ example, a physicist may be interested in what is the meaning of the "matter" 1-form fields in the equations. The answer is as follows. If you consider the conformal field theory, corresponding to open strings in dimension $2D$ on the torus, it appears to be a logarithmic (see e.g. \cite{pol}) VOA, and the Lian-Zuckerman construction does not work there.  In order to get rid of logarithms we considered another VOA on this space, corresponding to $\beta$-$\gamma$ systems, with the deformed BRST operator. However, even with this deformation we didn't recover the original open string theory, but the one, where half of dimensions is $T$-$dualized$.  
We expect also the relation of our considerations to the so-called "Double Field Theory" introduced by Hull and Zwiebach \cite{dft}, where the similar structures, like the Courant algebroid appear in the context of T-duality in String Field Theory.\\

\noindent {\bf 6.4. Nontrivial metric and a B-field. Einstein equations.} Another question one can ask, being motivated by the example with $\beta$-$\gamma$ system is about the meaning of the Yang-Mills $A_{\infty}$-algebra on a manifold with nontrivial metric and a $B$-field from the VOA point of view. So far, we have the construction only with a flat metric. We propose the following solution to this problem. We conjecture that there should be an 
$L_{\infty}$-algebra action on the LZ homotopy commutative algebra, in such a way that the flat background deformation can be treated as a deformation related to the Abelian $L_{\infty}$-subalgebra. From the form of $\eta$-deformation, one may suggest that this $L_{\infty}$-algebra comes from a tensor product of two LZ algebras. There is an indication of the existence of such $L_{\infty}$-algebra in \cite{lmz}, \cite{zeit2}, \cite{zeit3}. The Maurer-Cartan equation for this homotopy Lie algebra should be equivalent to the Einstein equations with external fields in the case of $\beta$-$\gamma$ VOA. We will address this question in \cite{zeitbeta}.

\section*{Acknowledgements}

I am indebted to I.B. Frenkel, M.M. Kapranov and M. Movshev for valuable discussions. 
I am very grateful to J. Stasheff for reading the draft version of this article and his comments. 
I would like to thank the referees for their valuable comments and remarks. 
I am grateful to the organizers of Simons Workshop'09, where this work was partly done.

\end{document}